\documentclass[12pt]{article}
\usepackage{latexsym}
\usepackage{amsfonts}
\usepackage{amsmath}
\usepackage{psfrag}
\usepackage{graphicx}
\usepackage{amssymb} 
\usepackage{bbold}
\usepackage{pstricks,pst-node,pst-tree}
\usepackage{array}
\usepackage{rotating}

\begin{document}
\centerline{\large{\bf A quantum Mermin--Wagner theorem}} 
\centerline{\large{\bf for quantum rotators on 2D graphs}}

\vspace{3mm}
\centerline{\bf Mark Kelbert}

\vspace{1mm}
\centerline{Department of Mathematics, Swansea University, UK}
\vspace{1mm}
\centerline{IME, University 
of S\~ao Paulo, Brazil}
\vspace{1mm}
\centerline{m.kelbert@swansea.ac.uk}
\vspace{3mm}
\centerline{\bf Yurii Suhov}
\vspace{1mm}
\centerline{StatsLab, DPMMS, University of Cambridge, UK}
\vspace{1mm}
\centerline{IME, University 
of S\~ao Paulo, Brazil}
\vspace{1mm}
\centerline{ITP, RAS, Moscow, Russia}
\vspace{1mm}
\centerline{yms@statslab.cam.ac.uk}
\vspace{3mm}

\def\ovr{\overline r}

\def\ttd{\tt d}
\def\ttg{\tt g}
\def\ttG{\tt G}
\def\bttg{\mbox{\boldmath${\ttg}$}}

\def\wt{\widetilde}
\def\wth{\wt h}
\def\wti{\wt i}
\def\wtj{\wt j}
\def\wtx{\wt x}
\def\wty{\wt y}

\def\wtalpha{\widetilde\alpha}
\def\gam{\gamma}
\def\Gam{\Gamma}
\def\bXi{\mbox{\boldmath${\Xi}$}}
\def\lam{\lambda}
\def\Lam{\Lambda}
\def\bmu{\mbox{\boldmath${\mu}$}}
\def\bnu{\mbox{\boldmath${\nu}$}}
\def\bphi{\mbox{\boldmath${\phi}$}}
\def\vphi{\varphi}
\def\bvphi{\mbox{\boldmath${\vphi}$}}
\def\bpsi{\mbox{\boldmath${\psi}$}}
\def\btau{\mbox{\boldmath${\tau}$}}
\def\om{\omega}
\def\Om{\Omega}
\def\bom{{\mbox{\boldmath$\om$}}}
\def\bOm{\mbox{\boldmath${\Om}$}}

\def\omu{\overline\mu}
\def\obmu{\overline\bmu}
\def\ovphi{\overline\vphi}
\def\obvphi{\overline\bvphi}
\def\obpsi{\overline\bpsi}

\def\oom{\overline\omega}
\def\oOm{\overline\Omega}
\def\sbom{\mbox{\sl${\bom}$}}
\def\obom{\overline\bom}
\def\0bom{{\bom}^0}
\def\0obom{{\obom}^0}
\def\nbom{{\bom}_n}
\def\0nbom{{\bom}_{n,0}}
\def\n*bom{{\bom}^*_{(n)}}

\def\utheta{\underline\theta}

\def\wtbom{\widetilde\bom}
\def\whbom{\widehat\bom}
\def\oom{\overline\om}
\def\wtom{\widetilde\om}

\def\obOm{\overline\bOm}
\def\whbOm{\widehat\bOm}
\def\wtbOm{\widetilde\bOm}

\def\oomega{\overline\omega}
\def\oUpsilon{\overline\Upsilon}
\def\wtomega{\widetilde\omega}
\def\wtheta{\widetilde\theta}

\def\fB{\mathfrak B}
\def\fG{\mathfrak G}
\def\fW{\mathfrak W}

\def\rd{\rm d}
\def\rr{\rm r}
\def\rx{\rm x}
\def\ry{\rm y}
\def\rtr{\rm{tr}}

\def\bcH{\mbox{\boldmath${\cH}$}}

\def\cl{\centerline} 

\def\cA{\mathcal A} 
\def\cB{\mathcal B}
\def\cC{\mathcal C}
\def\cE{\mathcal E} 
\def\cF{\mathcal F} 
\def\cH{\mathcal H} 
\def\cK{\mathcal K} 
\def\cL{\mathcal L} 
\def\cT{\mathcal T}
\def\cV{\mathcal V} 
\def\cW{\mathcal W}
\def\ocW{\overline\cW}

\def\bbC{\mathbb C}  
\def\bbE{\mathbb E}    
\def\bbP{\mathbb P}
\def\bbQ{\mathbb Q}
\def\bbR{\mathbb R}
\def\bbS{\mathbb S}
\def\bbT{\mathbb T}
\def\bbZ{\mathbb Z}

\def\obbP{\overline\bbP}

\def\oJ{\overline J}
\def\oP{\overline P}

\def\ba{\mathbf a}
\def\bA{\mathbf A}
\def\bB{\mathbf B}
\def\bF{\mathbf F}
\def\bH{\mathbf H}
\def\bI{\mathbf I}
\def\bN{\mathbf N}
\def\bR{\mathbf R}
\def\bU{\mathbf U}
\def\bX{\mathbf X}
\def\bx{\mathbf x}
\def\oB{\overline B}
\def\obx{\overline{\mathbf x}}
\def\ox{\overline x}
\def\wtbx{\widetilde\bx}
\def\uk{\underline k}
\def\un{\underline n}
\def\ux{\underline x}
\def\uy{\underline y}
\def\wtux{\widetilde\ux}
\def\uX{\underline X}

\def\bF{\mathbf F}
\def\bh{\mathbf h}
\def\by{\mathbf y}
\def\bn{\mathbf n}
\def\uy{\underline y}
\def\uY{\underline Y}
\def\uz{\underline z}
\def\bz{\mathbf z}
\def\uv{\underline v}
\def\uth0{{\underline \eta}_0}
\def\dist{\textrm{dist}}
\def\diy{\displaystyle}
\def\ov{\overline}

\def\wh{\widehat}

\def\wtz{\widetilde z}

\def\bI{\mathbf I}
\def\bK{\mathbf K}
\def\bP{\mathbf P}
\def\bV{\mathbf V}
\def\oW{\overline W}
\def\ofW{\overline\fW}
\def\bY{\mathbf Y}
\def\bfB{\mbox{\boldmath${\fB}$}}
\def\bfG{\mbox{\boldmath${\fG}$}}
\def\bfW{\mbox{\boldmath${\fW}$}}

{\bf Abstract.} This is the first of a series of papers considering
symmetry properties of quantum systems over 2D graphs or manifolds, 
with continuous spins, in the spirit of the Mermin--Wagner theorem
\cite{MW}. In the model considered here (quantum rotators)
the phase space of a single spin is a $d-$dimensional torus 
$M$, and spins (or particles) are attached to sites of a graph
$(\Gam ,\cE )$ satisfying a special bi-dimensionality property.
The kinetic energy part of the Hamiltonian is minus a half of 
the Laplace operator $-\Delta /2$ on $M$.
We assume that the interaction potential is C$^2$-smooth and
invariant under the action of a connected Lie group 
${\ttG}$ (i.e., a Euclidean space ${\bbR}^{d'}$ or 
a torus $M'$ of dimension $d'\leq d$) on $M$ preserving 
the flat Riemannian metric. A part of our approach is to give a
definition (and a construction) of a class of infinite-volume 
Gibbs states for the systems under consideration (the class $\fG$). 
This class contains the so-called limit Gibbs states, with or without
boundary conditions. We use
ideas and techniques originated from papers \cite{DS}, \cite{P}, 
\cite{FP}, \cite{SS} and \cite{ISV}, in combination
with the Feynman--Kac representation, to
prove that any state lying in the class $\fG$ (defined in the text) is 
${\ttG}$-invariant. An example is given where the interaction
potential is singular and there exists a Gibbs state which 
is not ${\ttG}$-invariant.

In the next paper under the same title we establish a similar result
for a bosonic model where particles can jump from a vertex $i\in\Gam$
to one of its neighbors (a generalized Hubbard model).
\vskip .5 truecm

\vskip .5 truecm

\vskip 2 truecm

\quad{\bf 1. Introduction. Existence and invariance of 
a limiting} 

\qquad\; {\bf Gibbs state}
\vskip 1 truecm

This work had been motivated, on the one 
hand, by a spectacular success on Mermin--Wagner
type theorems achieved in the past 
for a broad class of two-dimensional classical and quantum systems (see 
the bibliography quoted below) and,
on the other hand, by a recognised progress in experimental quantum physics 
creating and working with thin materials like graphene. 
The main dissatisfaction with published rigorous results in this area
stems for us from the fact that a
natural class of quantum models remained uncovered. These are systems
where the Hamiltonian contains a kinetic energy part given by a Laplacian. 
A serious problem here is that the finite-volume
Hamiltonians are unbounded operators. As a result, the construction of the
infinite-volume dynamical group encounters difficulties
(it works fine for simplified quantum spins models (like Heisenberg's)
where the phase space of a spin is finite-dimensional). Consequently, the
KMS-definition of an infinite-volume Gibbs state lacks substence for this
class of models, apart from the non-interacting case. (At least this
is the situation as we know it at the time of writing these lines.) A 
consistent definition of an infinite-volume Gibbs state is a cornerstone 
for the concept of a phase transition (as a non-uniqueness phenomenon);
it is precisely this concept that makes the Mermin--Wagner theorem 
important (and elegant). 
\vskip .5 truecm

{\bf 1.1. Bi-dimensional graphs.}
In the present paper we focus on Mermin\\ --Wagner type result for
a quantum bosonic system with continuous spins, 
over a denumerable graph $(\Gam,\cE)$ (with a 
vertex set $\Gam$ and an edge set $\cE\subset\Gam\times\Gam$).    
The graph will be assumed to satisfy
a specific bi-dimensional property generalising properties 
of `regular' lattices such as a square lattice $\bbZ^2$ or 
a triangular lattice $\bbZ^2_{\triangle}$. Cf. Eqns (1.1.1), (1.1.2) below. 
(Graphene is clearly
a regular 2D lattice; however, the whole theoretical methodology could 
be examined in the context of a more general graph with a distinct 
bi-dimensionality property.)
More precisely, we assume that $(\Gam ,\cE)$  
has the property that whenever edge $(j',j'')\in\cE$,
the reversed edge $(j'',j')\in\Upsilon$ as well. Furthermore, 
$(\Gam ,\cE)$ is without multiple edges and  
has a bounded degree. The latter means that the number of 
edges $(j,j')$ with 
a fixed initial or terminal vertex is uniformly bounded:
$$\begin{array}{r}
\sup\Big[\max\big(\sharp\,\{j'\in\Gam:\;(j,j')\in\Upsilon\},
\qquad\qquad{}\\
\sharp\,\{j'\in\Gam:\;(j',j)\in\cE\}\big):\;j\in\Gam\Big]
<\infty.\end{array}\eqno{(1.1.1)}$$
The bi-dimensionality 
property is expressed in the bound 
$$0<\sup \left[\frac{1}{n}\,\sharp\,\Sigma (j,n):\;
j\in\Gam ,\,n=1,2,\ldots\right]
<\infty \eqno (1.1.2)$$  
where $\Sigma (j,n)$ denotes the set of vertices in $\Gam$ at 
graph distance $n$ from site $j\in\Gam$ (a sphere of radius
$n$ about $j$):
$$\Sigma (j,n)=\{j'\in\Gam:\;{\tt d}(j,j')=n\}.\eqno (1.1.3)$$ 
(The graph distance ${\tt d}(j,j')={\tt d}_{\Gam,\cE}(j,j')$ 
between sites $j,j'\in\Gam$ is defined as the minimal length of a 
path on $(\Gam ,\cE )$ joining $j$ and $j'$.) This implies that
the cardinality of the ball
$$\Lam (j,n)=\{j'\in\Gam:\;{\tt d}(j,j')\leq n\}.\eqno (1.1.4)$$ 
grows at most quadratically in $n$.
\vskip .5 truecm

{\bf 1.2. The phase space and the group action.} We consider the following 
model. With each site (vertex)
$j\in\Gam$ there is associated a Hilbert space $\cH$ realized 
as $L_2(M,v)$ where $M$ is a compact Riemannian 
manifold; $v$ stands for the induced Riemannian volume. 
In this paper we assume that $M$ is a $d-$dimensional torus $\bbR^d/\bbZ^d$.
However, parts of the argument which can be easily done for a general manifold
are conducted without referring to the specific case of the torus. 
(The full generalization of 
the main results for a general compact Riemannian manifold will be
discussed elsewhere.)
Physically, $\cH$ is the phase space of a quantum spin `attached' to a single 
site of the graph and $M$ is its classical prototype. We assume 
that a connected Lie group ${\ttG}$ is given, acting on $M$ and 
preserving the flat metric on $M$. Transitivity of the action is not needed,
hence ${\ttG}$ is itself a torus or a Euclidean space of dimension $d'\leq d$. 
The action is generally referred to as 
$$({\ttg},x)\in{\ttG}\times M\mapsto{\ttg}x\in M.\eqno (1.2.1)$$
An alternative is the additive form of writing: we represent an element
${\ttg}\in{\ttG}$ with a $d$-dimensional vector 
$$\utheta =\theta A
$$
where $\theta\in\bbR^{d'}/\bbZ^{d'}$ is a vector of dimension $d'$
and $A$ is a $d'\times d$ matrix with rational coefficient of the column 
rank $d'$. The action is then written as
$$(\utheta ,x)\mapsto x+\utheta\;{\rm{mod}}\;1.\eqno (1.2.2)$$
We will use both forms: the multiplicative form (1.2.1) makes formulas 
shorter whereas the additive one is more convenient in technical calculations. 

A physical example of a system of the above type is a `frustrated'
2D crystal lattice. Here some `heavy' atoms or ions are placed at
the vertices of a graph, and each atom possesses a light bosonic
particle moving according to standard rules of Quantum Mechanics.
A more complicated model arises when the number of particles
is not fixed, and they can `jump' from one vertex to another; see \cite{KS2}.

Another example emerges from quantum gravity: cf. \cite{KSY1}, \cite{KSY2}. Here,
a graph is random and emerges from (random) triangulations of a 
$1+1$-dimensional space-time complex. (The paper 
\cite{KSY1} deals with classical spins on random triangulations;
a quantum version of the model is treated in \cite{KSY2}.) Classical models  on 
general graphs with a variable structure have been treated in a recent paper
\cite{KoKP}. 

If $\Lam$ is a finite subset in $\Gam$ then the phase space of the 
quantum system over $\Lam$ is $\cH_{\Lam}:=\cH^{\otimes\Lam}$, 
the Hilbert space $L_2(M^{\Lam},
v^{\Lam})$. Here and below the superscripts $\otimes\Lam$ and $\Lam$  
mean, respectively, the tensor product 
of copies of $\cH=L_2(M,v)$ and the Cartesian
products of copies of $M$ and $v$, labelled by sites $j\in\Lam$.
Formally, elements of $\cH^{\otimes\Lam}$ are (complex) functions
$$\bx_\Lam =(x(j),\,j\in\Lam )\in 
M^{\Lam}\mapsto\phi (\bx_\Lam )\in\bbC$$ 
considered modulo a set of $v^{\Lam}$-measure $0$, with the 
standard norm and the scalar product
$$\|\phi\|=\left(\diy\int_{M^{\Lam}}|\phi (\bx_\Lam )|^2
\prod_{j\in\Lam}v({\rd}x(j))\right)^{1/2},$$
and
$$\langle\phi_1 ,\phi_2\rangle =
\diy\int_{M^{\Lam}}\phi_1(\bx_\Lam )
\overline{\phi_2(\bx_\Lam )}\prod_{j\in\Lam}v({\rd}x(j)).$$
The argument $\bx_\Lam\in M^{\Lam}$ represents a classical
configuration of particles in $\Lam$.
Physically, this setting leads to a bosonic nature of the models under
consideration.  

The action of ${\ttG}$ determines unitary operators $U_\Lam ({\ttg})$,
${\ttg}\in{\ttG}$, in $\cH_\Lam$:
$$U_\Lam(\ttg)\phi (\bx_\Lam )=\phi ({\ttg}^{-1}\bx_\Lam)
\;\hbox{where}\;{\ttg}^{-1}\bx_\Lam=\{{\ttg}^{-1}x(j), j\in\Lambda\}.\eqno (1.2.3)$$
\vskip .5 truecm

{\bf 1.3. The Hamiltonian of the model and assupmtions about the potential.} 
A standard form of the kinetic energy operator for an individual 
spin is $-\Delta/2$ where $\Delta$ stands for the Laplacian 
operator in $\cH$. We also assume that a two-body interaction potential is given, 
which is described by a real-valued function 
$$\big((x',j'),(x'',j'')\big)\mapsto J({\ttd}(j',j''))V(x',x'').\eqno (1.3.1)$$ 
In the main body of the paper we assume that the (real) function 
$(x',x'')\in M\times M\mapsto V(x',x'')$ is 
of class C$^2$, although in one particular result, Theorem 1.4, we consider
an `opposite' situation of a singular potential. (In a 
forthcoming paper, we will address in detail the case of quantum models with 
non-smooth potentials.)  

More precisely, in Theorems 1.1, 1.2, 3.1--3.2 and Corollary 3.3 below
we assume that the function $V$ and its first and second derivatives 
$\nabla_{\rx}V$
and $\nabla_{{\rx}_1}\nabla_{{\rx}_2}V$ 
satisfy the uniform bounds: $x',x''\in M$
$$
|V(x',x'')|, \left|\nabla_{\rx}V(x',x'')\right|, 
\left|\nabla_{{\rx}_1}\nabla_{{\rx}_2}V(x',x'')\right|\leq {\ov V}.
\eqno{(1.3.2)}$$ 
Here ${\rx}$, ${\rx}_1$ and ${\rx}_2$ run through the arguments 
$x',x''\in M$; $|\;\;|$ stands for the absolute value of a real scalar or 
the norm of a real vector, and ${\ov V}\in (0,+\infty )$ is a constant.
Next, the function $J:\;r\in (0,\infty )\mapsto
J(r)\geq 0$ is assumed monotonically non-increasing with $r$ and obeying
the relation ${\oJ}(l)\to 0$ as $l\to\infty$ where  
$${\oJ}(l)=\sup\,\left[\sum\limits_{j''\in\Gam}J ({\tt d}(j',j'')){\mathbf 1}({\tt d}(j',j'')\geq l):\;
j'\in\Gam\right]<\infty .\eqno{(1.3.3)}$$
Additionally, let the interaction potential be such that
$$ J^*=\sup\left[\sum_{j''\in \Gam} J(d(j',j'')) d(j',j'')^2: j'\in\Gam\right]<\infty.
\eqno (1.3.4)$$
Next, we assume that the function $V$ 
is ${\ttg}$-invariant: 
$$V(x,x')=V({\ttg}x,{\ttg}x'),\;\;\forall\;x,x'\in M,\; 
{\ttg}\in {\ttG}.\eqno{(1.3.5)}$$

The Hamiltonian $H_\Lam$ of the system over a finite
set $\Lam\subset\Gam$ acts on functions $\phi\in 
\cH^{\otimes\Lam}$: given $\bx_{\Lam}=(x(j),\,j\in\Lam )
\in M^{\Lam}$,
$$\big(H_{\Lam}\phi\big)(\bx_{\Lam})=
\frac{1}{2}\left[-\sum_{j\in\Lam}\Delta_j+
\sum_{j,j'\in\Lam\times\Lam}J({\ttd}(j,j')) V(x_j,x_{j'})\right]
\phi (\bx_{\Lam}).\eqno{(1.3.6)}$$ 
Here $\Delta_j$ stands for the Laplace operator in 
variable $x(j)\in M$. A more general concept is a Hamiltonian $H_{\Lam |
\obx_{\Gam'\setminus\Lam}}$ in the external field generated by a
(finite or infinite) configuration $\obx_{\Gam'\setminus\Lam}
=\{\ox_{j'},\,j'\in\Gam'\setminus\Lam\}
\in M^{\Gam'\setminus\Lam}$ where $\Gam'\subseteq\Gam$ is
a (finite or infinite) collection of vertices. Namely,
$$\begin{array}{r}\diy\big(H_{\Lam |\obx_{\Gam'\setminus\Lam}}
\phi\big)(\bx_{\Lam})=
\bigg[-\frac{1}{2}\sum\limits_{j\in\Lam}\Delta_j+\frac{1}{2}
\sum\limits_{(j,j')\in\Lam\times\Lam}J({\ttd}(j,j'))V(x_j,x_{j'})\qquad{}\\
\diy +\sum\limits_{(j,j')\in\Lam\times (\Gam'\setminus\Lam)}J({\ttd}(j,j'))V(x_j,\ox_{j'})\bigg]
\phi (\bx_{\Lam}).\end{array}\eqno{(1.3.7)}$$ 

Summarizing, the model considered in this paper can be called a
system of quantum rotators on a bi-dimensional graph.
\vskip .5 truecm

{\bf 1.4. Properties of limiting Gibbs states.} Throughout the paper, 
we use a number of well-known facts (properties (i)--(iv) and 
(a)--(c) below) related to 
operators $H_{\Lam}$ and $H_{\Lam |\obx_{\Gam'\setminus\Lam}}$ 
which can be extracted, e.g., from Refs \cite{BR}, \cite {Gi}, 
\cite{RS}, \cite{S1}. (i) Under the above assumptions, 
operators $H_{\Lam}$ and $H_{\Lam |\obx_{\Gam'\setminus\Lam}}$
are self-adjoint
(on the natural domains) in $\cH^{\otimes\Lam}$,  bounded from 
below and have a discrete 
spectrum. (ii) Moreover, $\forall$ $\beta >0$, $H_{\Lam}$ and
$H_{\Lam |\obx_{\Gam'\setminus\Lam}}$ give rise to 
positive-definite trace-class operators $\exp\,\big[-\beta H_{\Lam}\big]$
and $\exp\,\left[-\beta H_{\Lam |\obx_{\Gam'\setminus\Lam}}\right]$.
(iii) In turn, this gives rise to Gibbs states 
$\vphi_\Lam =\vphi_{\beta ,\Lam}$ and $\vphi_{\Lam |\obx_{\Gam'\setminus\Lam}}=
\vphi_{\beta ,\Lam |\obx_{\Gam'\setminus\Lam}}$, 
at temperature $\beta^{-1}$ in volume 
$\Lam$. These are linear positive
normalized functionals on the C$^*$-algebra $\fB_\Lam$ of bounded operators
in space $\cH_\Lam$:
$$\vphi_\Lam (A)={\rtr}_{\cH_{\Lam}}\big(R_\Lam A\big),\;
\vphi_{\Lam |\obx_{\Gam'\setminus\Lam}} (A)
={\rtr}_{\cH_{\Lam}}\big(R_{\Lam |\obx_{\Gam'\setminus\Lam}}A\big),\;
\;\;A\in\fB_\Lam ,\eqno{(1.4.1)}$$
where 
$$R_\Lam =\frac{\exp\,\big[-\beta H_{\Lam}\big]}{\Xi_{\beta ,\Lam}}\;
\hbox{ with }\;\;
\Xi_{\beta ,\Lam}
={\rtr}_{{\mathcal H}_\Lam}\big(\exp\,\big[-\beta H_\Lam\big]\big)
\eqno (1.4.2)$$
and
$$\begin{array}{l}\diy R_{\Lam |\obx_{\Gam'\setminus\Lam}} =
\frac{\exp\,\big[-\beta H_{\Lam |\obx_{\Gam'\setminus\Lam}}\big]}{
\Xi_{\beta ,\Lam |\obx_{\Gam'\setminus\Lam}}}\\
\diy\qquad\hbox{with }\;\;
\Xi_{\beta ,\Lam |\obx_{\Gam'\setminus\Lam}}
={\rtr}_{{\mathcal H}_\Lam} 
\big(\exp\,\big[-\beta H_{\Lam |\obx_{\Gam'\setminus\Lam}}\big]\big).
\end{array}\eqno (1.4.3)$$

(iv) Let $\fB$ stand for the C$^*$-algebra of bounded operators in Hilbert space
$\cH$. For $\Lam^0\subset\Lam$, the
representations $\fB_\Lam =\fB^{\otimes\Lam}$ and
$\fB_{\Lam^0}=\fB^{\otimes\Lam^0}$ identify
$\fB_{\Lam^0}$ 
with the C$^*$-sub-algebra in $\fB_\Lam$ formed by the operators of
the form $A_0\otimes I_{\Lam\setminus\Lam^0}$ where
$I_{\Lam\setminus\Lam^0}$ is the unit operator in $\cH_{\Lam\setminus
\Lam^0}$. Accordingly, the restriction
$\vphi_\Lam^{\Lam^0}$ of state $\vphi_\Lam$ to 
C$^*$-algebra $\fB_{\Lam^0}$
is given by
$$\vphi_\Lam^{\Lam^0}(A_0)={\rtr}_{\cH_{\Lam^0}}
\big( R_\Lam^{\Lam^0}A_0\big),\;\;A_0\in\fB_{\Lam^0},\eqno{(1.4.4)}$$
where
$$R_\Lam^{\Lam^0}=
{\rtr}_{\cH_{\Lam\setminus\Lam^0}}R_\Lam.\eqno{(1.4.5)}$$
Clearly, operators $R_\Lam^{\Lam^0}$ are positive-definite and
have ${\rtr}_{\cH_{\Lam^0}}R_\Lam^{\Lam^0}=1$. They also
satisfy the compatibility property: $\forall$ $\Lam^0\subset\Lam^1
\subset\Lam$,
$$R_\Lam^{\Lam^0}=
{\rtr}_{\cH_{\Lam^1\setminus\Lam^0}}R_\Lam^{\Lam^1}.
\eqno{(1.4.6)}$$
Furthermore, in a similar fashion one can define functionals
$\vphi_{\Lam |\obx_{\Gam'\setminus\Lam}}^{\Lam^0}$
and operators 
$R_{\Lam |\obx_{\Gam'\setminus\Lam}}^{\Lam^0}$, with the same 
properties.

Below we denote by $\Lam\nearrow\Gam$ the net of finite subsets
of $\Gam$ ordered by inclusion. A convenient example of an increasing  
sequence in this net, eventually covering the entire $\Gam$, is 
formed by sets $\Lam (j,n)$, $n=1,2,\ldots$ (balls in the graph distance); 
see (1.1.4).

We prove in this paper the following results:
\medskip\medskip

{\bf Theorem 1.1.} {\sl For all given $\beta\in (0,\infty )$ and a finite
$\Lam^0\subset\Gam$, operators $R_\Lam^{\Lam^0}$ 
form a compact sequence in the trace-norm topology in $\cH_{\Lam^0}$
as $\Lam\nearrow\Gam$. Furthermore, given any family of 
(finite or infinite)
sets $\Gam' =\Gam'(\Lam )\subseteq\Gam$ and particle configurations
$\obx_{\Gam'\setminus\Lam}$,
operators $R_{\Lam |\obx_{\Gam'\setminus\Lam}}^{\Lam^0}$ also 
form a compact sequence in the trace-norm topology. 

Moreover, any limiting point, $R^{\Lam^0}$, for 
$\left\{R_{\Lam |\obx_{\Gam'\setminus\Lam}}^{\Lam^0}\right\}$ 
is a positive definite operator of trace one which possesses the following 
invariance property: $\forall$ ${\ttg}\in{\ttG}$,
$$U_{\Lam^0}({\ttg})^{-1}R^{\Lam^0}U_{\Lam^0}({\ttg})=
R^{\Lam^0}.\quad\lhd\eqno (1.4.7)$$
}
\medskip

By invoking the diagonal process, we get a family $\{R^{\Lam^0}\}$ 
of positive definite operators
$R^{\Lam^0}$ in $\cH_{\Lam^0}$ of trace one, where $\Lam^0$ 
runs over finite subsets of $\Gam$,
with the following properties. (a) $\exists$ an increasing sequence of 
finite sets $\Lam_{n_k}\subset\Gam$ such that 
$\operatornamewithlimits{\bigcup}\limits_{k}\Lam_{n_k}=\Gam$ and
a sequence of sets $\Gam'_{n_k}\subseteq\Gam$ and particle configurations
$\obx_{\Gam'_{n_k}\setminus\Lam_{n_k}}$ such that  
for all finite set $\Lam^0$ the convergence in the trace-norm
holds:
$$R^{\Lam^0}=\lim_{k\to\infty}R_{\Lam_{n_k}|\obx_{\Gam'_{n_k}
\setminus\Lam_{n_k}}}^{\Lam^0}.
\eqno{(1.4.8)}$$
(b) $\forall$ finite subsets $\Lam^0$, $\Lam^1$ of $\Gam$,
with $\Lam^0\subset\Lam^1$,
$$R^{\Lam^0}=
{\rtr}_{\cH_{\Lam^1\setminus\Lam^0}}R^{\Lam^1}.\eqno{(1.4.9)}$$
(c) Such a family defines a state $\vphi$ of (that is, 
a linear positive normalized functional on) 
the $C^*$-algebra quasilocal observables $\diy{\fB}_\Gam =
\left({\fB}^0_\Gam\right)^-$. Here $*$-algebra
$\diy{\fB}^0_\Gam$ is the inductive limit $\diy
{\operatornamewithlimits{\lim\,\rm{ind}}\limits_{n\to\infty}}\;
\fB_{\Lam_n}$ and superscript $^-$ in the notation 
$\left({\fB}^0_\Gam\right)^-$ stands for the norm completion. 
See \cite{BR}.

The definition of the above state $\vphi$ is that $\forall$
finite $\Lam^0\subset\Gam$, 
$$\vphi (A)={\rtr}_{\cH_{\Lam^0}}R^{\Lam^0}A.$$ 
Reflecting its construction, we call $\vphi$
a limiting Gibbs state; Theorem 1.1 asserts that the set $\fG^0$ of limiting 
Gibbs states is non-empty. A straightforward corollary is
\medskip

{\bf Theorem 1.2.} {\sl Any limiting Gibbs state $\vphi\in\fG^0$ has the
following invariance property: $\forall$ finite $\Lam^0\subset\Gam$
any $A\in\fB_{\Lam^0}$ and ${\ttg}\in{\ttG}$,} 
$$\vphi (A)=\vphi (U_{\Lam^0}({\ttg})^{-1}
AU_{\Lam^0}({\ttg})).\quad\lhd\eqno (1.4.10)$$
\medskip

The proof of Theorem 1.1 
is based on the following Lemma.
\vskip .5 truecm

{\textbf{Lemma 1.1.}} {\sl Let $\rho_n(x,y)$ be a sequence of kernels
defining positive-definite
operators $R_n$ of trace class and with trace $1$ in a Hilbert
space $L_2(M,\nu )$ where $\nu (M)<\infty$. Suppose there
exists the following limit, uniform in $x,y\in M$:
$$\lim_{n\to\infty}\rho_n(x,y)=\rho (x,y),\eqno (1.4.11)$$
which defines a positive-definite trace-class operator $R$
of trace $1$. Then
$$\lim_{n\to\infty}\|R_n-R\|_{\rm{tr}}=0\eqno (1.4.12)$$
where $\|A\|_{\rm{tr}}={\rm{tr}}\big(AA^*\big)^{1/2}$. $\quad\lhd$}
\vskip .5 truemm

Lemma 1.1 appeared for the first time in the short note \cite{Su1}. For the 
reader's convenience we give a complete proof in Section 4.3. 
\vskip .5 truemm

{\bf Remark 1.1.} As usually with  Mermin--Wagner type assertions, Theorem
1.2 does not address the issue of phase transitions, viz., uniqueness
of a limiting Gibbs state. A matter of principle here is to determine 
within what class of states $\fG\supseteq\fG^0$ the invariance 
property still holds true. 
Such a class is introduced in the next section; it is related to the
Feynman--Kac representation of operator $\exp\,\big[-\beta H_{\Lam}\big]$.
\vskip 1 truecm

Throughout the paper we adopt the following notational agreement: symbol 
$\lhd$ marks the end of 
a statement and symbol $\Box$ the end of a proof. 
\vfill\eject

{\bf 2. The Feynman--Kac formula and DLR equations}
\vskip 1 truecm

{\bf 2.1. The Feynman--Kac (FK) representation for the partition 
function.} 
In this section we follow the approach developed in \cite{Gi}; see also
\cite{AKKR}. Our
first observation is that, under the above assumptions,
operator $\exp\,\big[-\beta H_{\Lam}\big]$
acts as an integral operator in variables ${\bx}_{\Lam}
=(x(j),\,j\in\Lam )\in M^{\Lam}$ and ${\by}_{\Lam}=
(y(j),\,j\in\Lam )\in M^{\Lam}$:
$$
\Big(\exp\,\big[-\beta H_{\Lam}\big]\phi\Big)(\bx_\Lam )
=\int_{M^\Lam} \prod\limits_{j\in\Lam}v({\rd}y(j)) 
{K}_{\beta ,\Lam}(\bx_{\Lam},\by_{\Lam}) 
\phi (\by_{\Lam}).\eqno{(2.1.1)}$$ 
The integral kernel ${K}_{\beta ,\Lam}(\bx_{\Lam},\by_{\Lam})$
admits a  Feynman--Kac (FK) integral representation 
$${K}_{\beta ,\Lam}(\bx_{\Lam},\by_{\Lam}) 
=\diy\int\limits_{\oW^\beta_{\bx_\Lam ,\by_\Lam}}
{\bbP}^{\beta}_{{\bx}_{\Lam},{\by}_{\Lam}}
({\rd}\obom_\Lam)
\exp\,\big[-h^{\Lam}(\obom_\Lam)\big]\eqno{(2.1.2)}$$ 
explained below.

In Eqn (2.1.2), $\oW^\beta_{\bx_\Lam ,
\by_\Lam}$ stands for the Cartesian product 
$\operatornamewithlimits{\times}\limits_{j\in\Lam}\oW^\beta_{x(j),y(j)}$. 
Next, the Cartesian  
factor $\oW^\beta_{x(j),y(j)}$ represents the space of 
continuous paths
$\oom_j$ in $M$, of time-length $\beta$ and with the end-points 
$x(j)$ and $y(j)$:
$$\begin{array}{c}
\oom_j :\tau\in [0,\beta ]\mapsto \oom_j(\tau )\in M:\;\;\oom_j(\,\cdot\,)\;
\hbox{continuous,}\\
\oom_j(0)=x(j),\;\;\oom_j(\beta )=y(j),\;\;j\in\Lam .
\end{array}$$
Correspondingly, $\obom_\Lam =(\oom_j,\;j\in\Lam )\in
\oW^\beta_{\bx_\Lam ,\by_\Lam}$ 
is a collection of continuous paths 
$\oom_j\in \oW^\beta_{x(j) ,y(j)}$, $j\in\Lam$. We will say
that $\obom_\Lam$ is a path configuration over $\Lam$. Further, 
${\bbP}^{\beta}_{\bx_{\Lam},\by_{\Lam}}$ is 
the product-measure on $\oW^\beta_{\bx_\Lam ,\by_\Lam}$:
$${\bbP}^{\beta}_{\bx_{\Lam},\by_{\Lam}}
({\rd}\obom_\Lam)=
\operatornamewithlimits{\times}\limits_{j\in\Lam}{\bbP}^\beta _{x(j) ,y(j)}
({\rd}\oom_j)\eqno{(2.1.3)}$$ 
where ${\bbP}^\beta _{x(j) ,y(j)}({\rd}\oom_j)$
is the (non-normalised) Wiener measure on $\oW^\beta_{x(j) ,y(j)}$
(the Brownian bridge in $M$, of time-length $\beta$, with endpoints 
$x(j)$ and $y(j)$). The measure ${\bbP}^\beta _{x(j) ,y(j)}({\rd}\oom_j)$
is defined on a standard sigma-algebra of subsets of $\oW^\beta_{x(j) ,y(j)}$
generated by cylinders, and the measure ${\bbP}^{\beta}_{\bx_{\Lam},\by_{\Lam}}$
on the corresponding sigma-algebra of subsets in $\oW^\beta_{\bx_\Lam ,\by_\Lam}$.
In future we do not always explicitely refer to the sigma-algebras where measures
under consideration are defined: their specification follows that of the
underlying spaces. 
 
Finally, for a path configuration $\obom_\Lam
=\{\omega_j,\;j\in\Lam\}$ over $\Lam$, 
$$h^{\Lam}(\obom_\Lam)=\sum\limits_{(j,j')\in\Lam\times\Lam} 
h^{j,j'}(\oom_j,\oom_{j'})\eqno{(2.1.4)}$$ 
where $h^{j,j'}(\oom_j,\oom_{j'})$ represents an integral along 
trajectories $\oom_j$ and $\oom_{j'}$: 
$$h^{j,j'}(\oom_j,\oom_{j'})=J({\ttd}(j,j'))
\int\limits_0^{\beta}{\rd}\tau 
\,V\big(\oom_j(\tau),\oom_{j'}(\tau)\big).\eqno{(2.1.5)}$$ 
It is convenient to think that $h^{j,j'}(\oom_j,\oom_{j'})$ yields the
`energy of interaction' between trajectories $\oom_j$ and $\oom_{j'}$,   
and $h^{\Lam}(\obom_\Lam)$ equals the `full potential energy'
of the path configuration $\obom_{\Lam}$.

Furthermore, the trace ${\rtr}_{\cH_{\Lam}}\exp\;
\big[-\beta H_{\Lam}\big]$ (the partition function in $\Lam$) 
is finite and equals $\Xi_{\beta ,\Lam}$ where
$$\Xi_{\beta ,\Lam}
=\int_{M^\Lam}\prod\limits_{j\in\Lam}v({\rd}x(j))
{K}_{\beta,\Lam}(\bx_{\Lam},\bx_{\Lam})<+\infty .\eqno{(2.1.6)}$$ 
Consequently, operator $R_\Lam$ from (1.4.2) (often called
the density matrix (DM) in $\Lam$) is given by its integral 
kernel $F_{\beta,\Lam}(\bx_{\Lam},\by_{\Lam})$ (the DM kernel,
DMK for short):
$$F_{\beta ,\Lam}(\bx_{\Lam},\by_{\Lam})
=\frac{1}{\Xi_{ \beta ,\Lam}}
{K}_{\beta ,\Lam}(\bx_{\Lam},\by_{\Lam}).\eqno{(2.1.7)}$$ 
\vskip .5 truecm

{\bf 2.2. The FK representation for the RDMK in a finite volume.}
The operator $R_\Lam^{\Lam^0}$ from (1.4.8), (1.4.9) (referred
to as a reduced DM, briefly, RDM)
is determined by its own integral kernel 
$F^{\Lam^0}_{\beta,\Lam}(\bx_{\Lam^0},\by_{\Lam^0})$
(the RDM kernel, shortly, RDMK):
$$
F^{\Lam^0}_{\beta,\Lam}(\bx_{\Lam^0},\by_{\Lam^0})
=\frac{\Xi_{ \beta ,\Lam\setminus\Lam^0}
(\bx_{\Lam^0},\by_{\Lam^0})}{\Xi_{ \beta ,\Lam}},\;\;\bx_{\Lam^0},\by_{\Lam^0}\in 
M^{\Lam^0}.\eqno (2.2.1)$$
Here the quantity $\Xi_{ \beta ,\Lam\setminus\Lam^0}
(\bx_{\Lam^0},\by_{\Lam^0})$
in the numerator yields a `partial' partition function
corresponding to the partial trace $
{\rtr}_{\cH_{\Lam\setminus\Lam^0}}$ in (1.4.6):
$$\begin{array}{r}\Xi_{ \beta ,\Lam\setminus\Lam^0}
(\bx_{\Lam^0},\by_{\Lam^0})=
\int_{M^{\Lam\setminus\Lam^0}}\prod\limits_{j\in\Lam
\setminus\Lam^0}v({\rd}z_j)\qquad{}\\
\times K_{\beta}(\bx_{\Lam^0}\vee\bz_{\Lam\setminus\Lam^0},
\by_{\Lam^0}\vee\bz_{\Lam\setminus\Lam^0}
),\end{array}\eqno{(2.2.2)}$$ 
where symbol $\vee$ means concatenation of configurations
(this notation will be repeatedly used below). 

It is convenient to use a brief notation
${\rd}{\bx}_{\Lam}$ 
for the product of the Riemannian volumes 
$\operatornamewithlimits{\times}\limits_{j\in\Lam}v({\rd}x(j))$. 
We will also omit, where possible, the argument/index $\beta$ from the 
notation (viz., by writing $\Xi_{\Lam}$ instead of $\Xi_{\beta ,\Lam}$). 
The above representations (2.1.1)--(2.1.7) allow us to associate with 
Gibbs state $\vphi_\Lam$ a probability distribution 
$\mu_{\Lam}$ on the set
$$W_{\Lam}=
\operatornamewithlimits{\bigcup}\limits_{\bx_\Lam\in M^{\Lam}}
\oW_{\bx_\Lam ,\bx_\Lam}\;\hbox{ where }\;
\oW_{\bx_\Lam ,\bx_\Lam}=\operatornamewithlimits{\bigcup}\limits_{j\in \Lambda}
\oW_{x(j),x(j)}.\;\eqno (2.2.3)$$  
the definition of this probability distribution is provided in forthcoming
paragraphs.

Pictorially, $W_{\Lam}$ is the space 
of collections of closed trajectories (loops) in $M$ issued from
and returning to (coinciding) specified endpoints; each loop
being assigned to a site $j\in\Lam$. 
Any such
loop collection can be written as a pair $(\bx_\Lam,\bom_\Lam )$.
Here $\bom_\Lam 
=\{\omega_j,\;j\in\Lam\}$ is a collection of loops
$\tau\in [0,\beta ]\mapsto\omega_j(\tau )$, where 
$\omega_j(0)=\omega_j(\beta )=x(j)$; a pair 
$(x(j),\om_j)$ is associated with site $j\in\Lam$. We will say that 
$(\bx_\Lam,\bom_\Lam )$ (and $\bom_\Lam$ when the reference to 
$\bx_\Lam$ is clear from the context) is a loop configuration over $\Lam$. Note 
the absence of the bar in this notation, stressing that 
$\bom_{\Lam}\in W_\Lam$ is a loop configuration as opposite to
a general path configuration 
$\obom_{\Lam}\in\oW_{\bx_\Lam ,\by_\Lam}\subset 
\oW_\Lam$ (again associated with sites 
$j\in\Lam$ (see Eqn (2.2.4) below)). 
More precisely, when appropriate, 
we will omit the bar in the notation $\oW_{x,y}$ and 
$\oW_{\bx_\Lam ,\by_\Lam}$ for $x=y$ or $\bx_\Lam =\by_\Lam$:
$$\oW_{x,x}=W_{x,x}\;\hbox{ and }\;
\oW_{\bx_\Lam ,\bx_\Lam}=W_{\bx_\Lam ,\bx_\Lam}.$$
Recall, we refer to $\bom_\Lam$ as a loop configuration 
and $\obom_{\Lam}$ as a path configuration in $\Lam$.
Next, we set:
$$\oW=
\operatornamewithlimits{\bigcup}\limits_{x,y\in M}\oW_{x,y}\;\hbox{ and }
\oW_{\Lam}=
{\operatornamewithlimits{\bigcup}\limits_{\bx_\Lam,
\by_\Lam\in M^{\Lam}}}
\oW_{\bx_\Lam,\by_\Lam}.\eqno{(2.2.4)}$$ 
\vskip .5 truecm

{\bf 2.3. The FK-DLR equations in a finite volume.}
The aforementioned probability distribution $\mu_{\Lam}$, on space 
$W_\Lam$, is absolutely 
continuous relative to the underlying product-measure $\nu_{\Lam}
(=\nu_{\beta ,\Lam})$, where 
$$\begin{array}{ll}{\rd}\nu_{\Lam}({\bx}_{\Lam},
\bom_{\Lam})&\diy =
{\operatornamewithlimits{\times}\limits_{j\in\Lam}}\left(v({\rd}x(j))
\times{\bbP}_{x(j),x(j)}({\rd}\om_j)\right)\\
\;&\diy ={\rd}\bx_\Lam\times
{\bbP}_{\bx_\Lam,\bx_\Lam}({\rd}\bom_{\Lam}).\end{array}\eqno{(2.3.1)}$$  
Here the measure ${\bbP}_{x(j),x(j)}({\rd}\om_j)$ is defined as a Brownian
bridge on manifold $M$ with the starting and end point $x(j)$.
Next, the Radon--Nikodym derivative (probability density function)
$\diy p_{\Lam}
({\bx}_{\Lam},\bom_\Lam):=
\frac{{\rd}\mu_{\Lam}({\bx}_{\Lam},
\bom_\Lam)}{{\rd}\nu_{\Lam}({\bx}_{\Lam},
\bom_{\Lam})}$ is of the form 
$$p_{\Lam}({\bx}_{\Lam},\bom_\Lam)
=\frac{1}{\Xi_{\Lam}}\exp\,\big[-h^{\Lam}(\bom_\Lam 
)\big]\eqno{(2.3.2)}$$  
where functional $h^{\Lam}(\bom_\Lam )$ 
has been defined in (2.1.4)--(2.1.5). It is convenient to treat $\mu_\Lam$ as 
a Gibbs probability measure for a `classical' spin system 
where `spins' are represented by loops affiliated with sites $j\in\Lam$. 

To shorten the notation we will omit henceforce the argument $\bx_\Lam$ 
and similar arguments from symbols like $p_{\Lam}
(\bx_{\Lam},\bom_\Lam)$, ${\rd}\nu_{\Lam}({\bx}_{\Lam},
\bom_\Lam)$ and ${\rd}\mu_{\Lam}({\bx}_{\Lam},
\bom_\Lam)$, bearing in mind that the initial/end-point configuration
$\bx_\Lam$ can be reconstructed from the loop configuration $\bom_\Lam$.  

Measure $\mu_{\Lam}$ defines 
a random field over $\Lam$ with realizations 
$\bom_\Lam =\{\om_j,\,j\in\Lam\}\in W_{\Lam}$ 
and has the following properties (I), (II).  
\vskip .1 truecm

(I) $\mu_{\Lam}$ {\sl satisfies the DLR equation over $\Lam$; 
cf. Eqn} (2.3.4) below. (Recall, $\Lam\subset\Gam$ is a finite set.) This 
means the following. Given $\Lam^0\subset\Lam$, 
let us agree to write $\bom^0$ for the loop configuration 
$\bom_{\Lam^0}\in W_{\Lam^0}$. 
Consider the partially integrated probability density
$$\begin{array}{r}\diy p_{\Lam}^{\Lam^0}(\bom^0)\diy :=
\int_{W_{\Lam\setminus\Lam^0}}
{\rd}\nu_{\Lam\setminus\Lam^0}(\bom_{\Lam\setminus\Lam^0})
p_{\Lam}(\bom^0\vee\bom_{\Lam\setminus\Lam^0})
\end{array}\eqno{(2.3.3)}$$ 
where $\bom^0\vee\bom_{\Lam\setminus\Lam^0}$ stands
for the concatenation of the two loop configurations yielding a loop
configuration over the whole of $\Lam$. Cf. Eqn (2.2.2).

Then, $\forall$ set $\Lam'$ such that $\Lam^0\subset\Lam'\subset
\Lam$, the density $p_{\Lam}^{\Lam^0}(\bom^0)$ obeys
$$\begin{array}{r}\diy p_{\Lam}^{\Lam^0}(\bom^0)=
\int_{W_{\Lam\setminus\Lam'}}
{\rd}\nu_{\Lam\setminus\Lam'}
(\bom_{\Lam\setminus\Lam'})\qquad\qquad\qquad{}\\
\diy\times
p_{\Lam}^{\Lam\setminus\Lam'}
(\bom_{\Lam\setminus\Lam'})
\frac{\Xi_{\Lam'\setminus\Lam^0}(\bom^0,
\bom_{\Lam\setminus\Lam'})}{
\Xi_{\Lam'}(\bom_{\Lam\setminus\Lam'})}.\end{array}\eqno{(2.3.4)}$$  
Here $p_{\Lam}^{\Lam\setminus\Lam'}
(\bom_{\Lam\setminus\Lam'})$
is the partially integrated density similar to (2.3.3):
$$\begin{array}{r}\diy p_{\Lam}^{\Lam\setminus\Lam'}
(\bom_{\Lam\setminus\Lam'}):=\int_{W_{\Lam'}}
{\rd}\nu_{\Lam'}(\wtbom_{\Lam'})
p_{\Lam}(\wtbom_{\Lam'}\vee\bom_{\Lam\setminus\Lam'}).
\end{array}\eqno{(2.3.5)}$$  
Further, the quantities $\Xi_{\Lam'\setminus\Lam^0}(\bom^0,
\bom_{\Lam\setminus\Lam'})$ and $\Xi_{\Lam'}
(\bom_{\Lam\setminus\Lam'})$ are given by the following integrals
$$\begin{array}{r}\diy
\Xi_{\Lam'\setminus\Lam^0}(\bom^0,\bom_{\Lam\setminus\Lam'})=
\int\limits_{W_{\Lam'\setminus\Lam^0}}
{\rd}\nu_{\Lam'\setminus\Lam^0}(\bom_{\Lam'\setminus\Lam^0})
\qquad\qquad{}\\
\times\exp\big[
-h^{\Lam'}(\bom^0\vee\bom_{\Lam'
\setminus\Lam^0}|
\bom_{\Lam\setminus\Lam'})\big]\end{array}\eqno{(2.3.6)}$$%
and
$$\begin{array}{r}\diy
\Xi_{\Lam'}(\bom_{\Lam\setminus\Lam'})=
\int\limits_{W_{\Lam'}}
{\rd}\nu_{\Lam'}(\bom_{\Lam'})\exp\big[
-h^{\Lam'}(\bom_{\Lam'}|
\bom_{\Lam\setminus\Lam'})\big].\end{array}\eqno{(2.3.7)}$$%
Next, given loop configurations $\bom_{\Lam'}=\{\om_j,j\in\Lam'\}\in
W_{\Lam'}$
and $\bom_{\Lam\setminus\Lam'}=\{\om_j,j\in\Lam\setminus\Lam'\}
\in M^{\Lam\setminus\Lam'}$, the functional
$h^{\Lam'}(\bom_{\Lam'}|
\bom_{\Lam\setminus\Lam'})$ in (2.3.7) is determined by
$$h^{\Lam'}
(\bom_{\Lam'}|\bom_{\Lam\setminus\Lam'})
=h^{\Lam'}(\bom_{\Lam'})
+h
(\bom_{\Lam'}|\,|\bom_{\Lam\setminus\Lam'})
\eqno{(2.3.8)}$$  
where the summand $h^{\Lam'}(\om_{\Lam'})$ is defined as in (2.1.4)
and $h
(\bom_{\Lam'}|\,|\bom_{\Lam\setminus\Lam'})$ is given by
$$h^{\Lam'}
(\bom_{\Lam'}|\,|\bom_{\Lam\setminus\Lam'})
=\sum\limits_{(j,j')\in\Lam'\times
\Lam\setminus\Lam'}h^{j,j'}(\om_j,\om_{j'})
\eqno{(2.3.9)}$$
with $h^{j,j'}(\om_j,\om_{j'})$ as in (2.1.5).

Finally, the functional
$h^{\Lam'}(\bom^0\vee\bom_{\Lam'
\setminus\Lam^0}|
\bom_{\Lam\setminus\Lam'})$ figuring in (2.3.6), for  
$\bom^0\in W_{\Lam^0}$, $\bom_{\Lam'}\in W_{\Lam'}$ and 
$\bom_{\Lam\setminus\Lam'}\in W_{\Lam\setminus\Lam'}$, is 
defined by similar formulas. We say that 
$h^{\Lam'}(\bom_{\Lam'}|
\bom_{\Lam\setminus\Lam'})$ and 
$h^{\Lam'}(\bom^0\vee\bom_{\Lam'
\setminus\Lam^0}|
\bom_{\Lam\setminus\Lam'})$ give the values of a `potential
energy' of the loop configurations $\bom_{\Lam'}$ and  
$\bom^0\vee\bom_{\Lam'}$ in the external field generated by 
$\bom_{\Lam\setminus\Lam'}$. 
In this context, $\Xi_{\Lam'}(\bom_{\Lam\setminus\Lam'})$ gives 
the partition function for loop configurations over the `intermediate volume' 
$\Lam'$ in an external potential field generated by the boundary 
condition $\bom_{\Lam\setminus\Lam'}\in W_{\Lam\setminus\Lam'}$.
Similarly, $\Xi_{\Lam'\setminus\Lam^0}
(\bom^0,\bom_{\Lam\setminus\Lam'})$ can be considered as 
the partition function in the `layer' 
$\Lam'\setminus\Lam^0$, with an `external' boundary condition 
$\bom_{\Lam\setminus\Lam'}\in W_{\Lam\setminus\Lam'}$ and an 
'internal' loop configuration $\bom^0\in W_{\Lam^0}$ (note that
$\bom^0$  enters the integral $\Xi_{\Lam'\setminus\Lam^0}
(\bom^0,\bom_{\Lam\setminus\Lam'})$
with its energy $h^{\Lam^0}(\bom^0)$). A straightforward fact is
that 
$$\int_{W_{\Lam^0}}{\rd}\nu_{\Lam^0}(\omega^0)
\Xi_{\Lam'\setminus\Lam^0}
(\bom^0,\bom_{\Lam\setminus\Lam'})=
\Xi_{\Lam'}(\bom_{\Lam\setminus\Lam'}).$$

In probabilistic terms, the DLR equation (2.3.4) is equivalent to the
following property. Consider the conditional distribution
${\rd}\mu_{\Lam}^{\Lam^0|\Lam\setminus\Lam'}(\om^0|\,
\bom_{\Lam\setminus\Lam'})$
induced by the probability measure $\mu_{\Lam}$, for a loop 
configuration $\bom^0$ over $\Lam^0$, conditioned by a loop configuration
$\bom_{\Lam\setminus\Lam'}$ over $\Lam\setminus\Lam'$. It 
is determined by the conditional probability density $\diy
p_{\Lam}^{\Lam^0\,|\,\Lam\setminus\Lam'}
(\bom^0\,|\bom_{\Lam\setminus\Lam'})
\diy :=\frac{{\rd}\mu_{\Lam}^{\Lam^0|\Lam\setminus\Lam'}
(\bom^0\,|\,\bom_{\Lam\setminus\Lam'})}{
{\rd}\nu_{\Lam^0}(\bom^0 )}$. The equivalent form of the DLR 
property means that this density has the form
$$\begin{array}{r}p_{\Lam}^{\Lam^0\,|\,\Lam\setminus\Lam'}
(\bom^0\,|\bom_{\Lam\setminus\Lam'})
\diy =\frac{\Xi_{\Lam'\setminus\Lam^0}
(\bom^0,\bom_{\Lam\setminus\Lam'})}{
\Xi_{\Lam'}(\bom_{\Lam\setminus\Lam'})}.
\end{array}\eqno{(2.3.10)}$$ 
In fact, $\mu_{\Lam}$ is the only measure that satisfies
the equations (2.3.5), (2.3.10). The name DLR (Dobrushin--Lanford--Ruelle) is
widely used in the classical statistical mechanics; see, e.g., \cite{G}.
\vskip .1 truecm

(II) {\sl The measure $\mu_\Lam$ determines the RDMK 
$F^{\Lam^0}_\Lam({\bx}^0,{\by}^0)$.}
Given $\Lam^0\subset\Lam$
and particle configurations ${\bx}^0,{\by}^0\in M^{\Lam^0}$, the 
RDMK $F^{\Lam^0}_{\Lam}({\bx}^0,{\by}^0)$ 
is defined by
$$F^{\Lam^0}_{\Lam}({\bx}^0,{\by}^0) =
\int\limits_{\oW_{\bx^0,\by^0}}
{\bbP}_{{\bx}^0,{\by}^0}({\rd}{\obom}^0){q}_{\Lam}^{\Lam^0}(\obom^0).
\eqno{(2.3.11)}$$ 
In turn, the functional ${q}_{\Lam}^{\Lam^0}(\obom^0)$ is determined
by the formula analogous to (2.3.4): $\forall$ $\Lam'$ such that 
$\Lam^0\subset\Lam'\subset\Lam$,
$$\begin{array}{r}\diy 
{q}_{\Lam}^{\Lam^0}(\obom^0)=
\int_{W_{\Lam\setminus\Lam'}}{\rd}\nu_{\Lam\setminus\Lam'}
(\bom_{\Lam\setminus\Lam'})\qquad\qquad{}\\
\diy \times 
p_\Lam^{\Lam\setminus\Lam'}(\bom_{\Lam\setminus\Lam'})
\,\frac{\Xi_{\Lam'\setminus\Lam^0}(\obom^0,
\bom_{\Lam\setminus\Lam'})}{
\Xi_{\Lam'}(\bom_{\Lam\setminus\Lam'})}\end{array}\eqno{(2.3.12)}$$ 
with quantities $p_\Lam^{\Lam\setminus\Lam'}(\bom_{\Lam
\setminus\Lam'})$,  $\Xi_{\Lam'\setminus\Lam^0}(\obom^0,
\bom_{\Lam\setminus\Lam'})$ and $\Xi_{\Lam^0}(
{\obom}_{\Lam\setminus\Lam^0})$ defined as in Eqns (2.3.5)--(2.3.8).
(The only difference is that a loop configuration $\bom^0$ in (2.3.7) has been
replaced with a more general path configuration $\obom^0$.)
Kernel $F^{\Lam^0}_{\Lam}({\bx}^0,{\by}^0)$ is often referred to as
the reduced DM kernel, briefly RDMK
(more precisely, the kernel of the DM in volume 
$\Lam$, reduced to $\Lam^0$). Accordingly, the functional 
${q}_{\Lam}^{\Lam^0}(\obom^0)$ can be called a 
reduced DM functional, briefly RDMF, for $\obom^0\in\oW_{\bx^0,\by^0}$. 
Similarly to (2.3.4), it will be convenient to write 
$$\begin{array}{r}{q}_{\Lam}^{\Lam^0}(\obom^0)=\diy
\int\limits_{W_{\Lam\setminus\Lam'}}{\rd}\nu_{\Lam\setminus\Lam'}
(\bom_{\Lam\setminus\Lam'})\qquad\qquad{}\\
\times p_\Lam^{\Lam\setminus\Lam'}(\bom_{\Lam\setminus\Lam'})
{q}_{\Lam}^{\Lam^0|\Lam\setminus\Lam'}(\obom^0|
\bom_{\Lam\setminus\Lam'})\end{array}\eqno{(2.3.13)}$$
where
$${q}_{\Lam}^{\Lam^0|\Lam\setminus\Lam'}(\obom^0|
\bom_{\Lam\setminus\Lam'})=
\frac{\Xi_{\Lam'\setminus\Lam^0}(\obom^0,
\bom_{\Lam\setminus\Lam'})}{
\Xi_{\Lam'}(\bom_{\Lam\setminus\Lam'})}.\eqno{(2.3.14)}$$
In analogy with Eqn (2.3.10), quantity ${q}_{\Lam}^{\Lam^0|
\Lam\setminus\Lam'}(\obom^0|\bom_{\Lam\setminus\Lam'})$ in (2.3.14)
can be called a conditional RDMF for a path configuration $\obom^0\in 
\oW_{\bx^0,\by^0}$, given a loop configuration 
$\bom_{\Lam\setminus\Lam'}\in W_{\Lam\setminus\Lam'}$.

Note that the RDMF $q^{\Lam^0}_\Lam(\obom^0)$ from (2.3.12), (2.3.13) satisfies the 
invariance relation 
$$
q^{\Lam^0}_\Lam(\obom^0)=q^{\Lam^0}_\Lam({\ttg}\obom^0),
\;\;{\ttg}\in{\ttG},\;\obom^0\in \oW_{\bx^0,\by^0}\eqno (2.3.15)
$$
where
$${\ttg}\obom^0=\{{\ttg}\oom_j,\;j\in\Lam^0\},\;\hbox{with}\;
({\ttg}\om_j)(\tau)={\ttg}(\om_j(\tau)),\,0\leq\tau\leq\beta .
\eqno (2.3.16)$$ 
Consequently, for RDMK $F^{\Lam^0}_{\Lam}({\bx}^0,{\by}^0)$  (see 
(2.3.11)), we have that 
$$
F^{\Lam^0}_{\Lam}({\bx}^0,{\by}^0) =
F^{\Lam^0}_{\Lam}({\ttg}{\bx}^0,{\ttg}{\by}^0),
\;\;{\ttg}\in{\ttG},\;\bx^0,\by^0\in M^{\Lam^0}\eqno (2.3.17)
$$
and for the RDM $R_\Lam^{\Lam^0}$
$$R_\Lam^{\Lam^0}=U_{\Lam^0}({\ttg})^{-1}R_\Lam^{\Lam^0}
U_{\Lam^0}({\ttg}),\;\;{\ttg}\in{\ttG}.\eqno (2.3.18)
$$

However, we will need to consider a similar RDMF 
$q^{\Lam^0}_{\Lam |\obx_{\Gam'\setminus\Lam}}(\obom^0)$
defined via operator $\exp\,\big[-\beta H_{\Lam |
\obx_{\Gam'\setminus\Lam}}\big]$ instead of 
$\exp\,\big[-\beta H_\Lam\big]$. (It can be called a conditional
RDMF, with the boundary condition $\obx_{\Gam'\setminus\Lam}$.) 
The invariance equation 
$$q^{\Lam^0}_{\Lam |\obx_{\Gam'\setminus\Lam}}(\obom^0)
=q^{\Lam^0}_{\Lam |\obx_{\Gam'\setminus\Lam}}({\ttg}\obom^0)
\eqno (2.3.19)$$
(i.e., an analog of (2.3.15)) fails; this makes the statements of
Theorems 1.1 and 1.2 non-trivial. (Of course, the covariance
property 
$$q^{\Lam^0}_{\Lam |\obx_{\Gam'\setminus\Lam}}(\obom^0)
=q^{\Lam^0}_{\Lam |{\ttg}\obx_{\Gam'\setminus\Lam}}({\ttg}\obom^0)$$
holds true but is useless for our purpose.)
\vfill\eject

{\bf 3. The class $\fG$ of Gibbs states in the infinite volume}
\vskip 1 truecm

{\bf 3.1. Definition of the class $\fG$.} The aim of this section is to define 
the invariance property
(2.3.17) (and consequently, property (2.3.18)) for 
functionals ${q}_{\Gam}^{\Lam^0}(\omega^0)$ (and related objects 
$F^{\Lam^0}_\Gam(\bx^0,\by^0)$ and $R_\Gam^{\Lam^0}$)
for the system in
an `infinite volume' (i.e., over the whole graph $(\Gam,\cE)$). That 
is, we want to prove that functional ${q}_{\Gam}^{\Lam^0}(\omega^0)$,
which we call inifinite-volume RDMF, obeys
$${q}_{\Gam}^{\Lam^0}(\obom^0)=
{q}_{\Gam}^{\Lam^0}({\ttg}\obom^0),\;{\ttg}\in{\ttG} ,\;\obom^0\in
\oW_{\Lam^0}.\eqno (3.1.1)$$  
The formal definition of infinite-volume RDMF 
${q}_{\Gam}^{\Lam^0}(\omega^0)$
(${q}^{\Lam^0}(\omega^0)$ for short)
related to the system over $(\Gam,\cE)$ requires additional 
constructions and will lead us to the definition of the aforementioned class
of states $\fG$; see below. 
At this point we state that the key step 
is to establish an asymptotical form of (2.3.19) for infinite-volume
conditional RDMF ${q}^{\Lam^0\,|\,\Gam\setminus\Lam'}
(\obom^0|\bom_{\Gam\setminus\Lam'})$ 
when set $\Lam'$ is `large enough'. 
In essense, we will prove that,
$\forall$ finite set $\Lam^0\subset\Gam$, 
$$\lim_{n\to\infty}\frac{{q}^{\Lam^0\,|\,\Gam\setminus\Lam (n)}
({\ttg}\obom^0|\bom_{\Gam\setminus\Lam (n)})}{
{q}^{\Lam^0\,|\,\Gam\setminus\Lam (n)}(\obom^0|
\bom_{\Gam\setminus\Lam (n)})}=1,\eqno (3.1.2)$$
uniformly in: (i) group element ${\ttg}\in{\ttG}$, (ii) path configuration 
$\obom\in\oW_{\Lam^0}$, (iii) an (infinite) loop configuration 
$\bom_{\Gam\setminus\Lam (n)}\in W_{\Gam\setminus\Lam (n)}$ 
representing an infinite-volume external boundary condition. Here and below  
$\Lam (n)=\Lam (o,n)$ and $\Sigma (n)$ mean the ball and the sphere of 
radius $n$ (cf. (1.1.4))
around a reference point $o\in\Gam$ (the choice of point $o$ will
not matter).

In fact, functional ${q}^{\Lam^0\,|\,\Gam\setminus\Lam (n)}
(\obom^0|\bom_{\Gam\setminus\Lam (n)})$ is itself defined as the limit
$${q}^{\Lam^0\,|\,\Gam\setminus\Lam (n)}(\obom^0|
\bom_{\Gam\setminus\Lam (n)})=\lim_{r\to\infty}
{q}_{\Lam (r)}^{\Lam^0\,|\,\Lam (r)\setminus\Lam (n)}(\obom^0|
\bom_{\Lam (r)\setminus\Lam (n)})\eqno (3.1.3)$$
where, for $\Lam^0\subset\Lam (n)\subset\Lam(r)$, the value
${q}_{\Lam (r)}^{\Lam^0\,|\,\Lam(r)\setminus\Lam (n)}(\obom^0|
\bom_{\Lam (r)\setminus\Lam (n)})$ has been determined in Eqn (2.3.14).

At this point it is appropriate to establish some probabilistic 
background. The Borel  sigma-algebra of subsets of the loop space 
$W$ is denoted by ${\fW}$. Given a finite subset $\Lam\subset\Gam$,
we obtain the induced sigma-algebra of subsets of $W_{\Lam}$ 
which is denoted by ${\fW}_{\Lam}$. Similarly, for a trajectory 
space $\oW$ the sigma-algebra $\ofW$ is defined which leads to 
the sigma-algebra $\ofW_{\Lam^0}$ of subsets 
in $\oW_{\Lam^0}$. For $\Lam'\subset\Lam$, the sigma-algebra 
${\fW}_{\Lam'}$ is naturally identified with
a sub-sigma-algebra of ${\fW}_{\Lam}$ which is denoted by the same symbol
${\fW}_{\Lam'}$). For the whole graph $\Gam$, we can introduce
the Cartesian product $W_\Gam$ considered as the countable set of  loop
configurations $\{\omega_j,\;j\in\Gam\}$, 
$\omega_j\in M$; as earlier, a loop $\omega_j$ is associated with site 
$j\in\Gam$. For a finite set $\Lam\subset\Gam$, the sigma-algebra 
${\fW}_{\Lam}$ can again be identified with the sigma-algebra
of subsets of $W_{\Gam}$; as before, it is convenient to use the same 
notation for both. The sigma-algebra ${\fW}_{\Gam}$ is defined as the 
smallest sigma-algebra of subsets of $W_{\Gam}$
containing ${\fW}_{\Lam}$ $\forall$ finite $\Lam\subset\Gam$.
In a similar fashion we define the sigma-algebra 
${\fW}_{\Gam\setminus\Lam^0}\subset
{\fW}_{\Gam}$ for a given (finite) set $\Lam^0\subset\Gam$;
as before, it is naturally identified with the sigma-algebra of subsets in 
$W_{\Gam\setminus\Lam^0}$.

Let us now define the class $\fG$ of states of the C$^*$-algebra $\fB$.
As before, a state $\vphi$ of $\fB$ is identified with a family of 
RDMs $R^{\Lam^0}=R^{\Lam^0}_\vphi$ 
where $\Lam^0$ is an arbitrary finite subset of $\Gam$;
each $R^{\Lam^0}$ is a positive definite operator in $\cH_{\Lam^0}$
of trace one, and the compatibility relation (1.4.9) holds true. In short,
for a state $\vphi\in\fG$ the RDMs $R^{\Lam^0}$ are integral operators
(see (3.1.4)),
with integral kernels $F^{\Lam^0}({\bx}^0,{\by}^0)$ satisfying 
(3.1.5)--(3.1.12), where the probability measure $\mu_\Gam$ obeys 
(3.1.13)--(3.1.17). Properties (3.1.5)--(3.1.17) are direct analogs of
the corresponding properties of the RDMKs $F^{\Lam^0}_\Lam({\bx}^0,{\by}^0)$ 
and $F^{\Lam^0}_{\Lam |\bx_{\Gam'\setminus\Lam}}({\bx}^0,{\by}^0)$ in
a finite volume $\Lam$.  
\vskip .5 truecm

Passing to the formal presentation, the RDM $R^{\Lam^0}$ is determined 
by its integral kernel $F^{\Lam^0}({\bx}^0,{\by}^0)$:
$$\big(R^{\Lam^0}\phi\big)({\bx}^0)=
\int_{M^{\Lam^0}}{\rd}\by^0 
F^{\Lam^0}({\bx}^0,{\by}^0) 
\phi ({\by}^0),\;\bx^0\in M^{\Lam^0}.\eqno (3.1.4)$$
In turn, the RDMK $F^{\Lam^0}({\bx}^0,{\by}^0)
=f_\vphi^{\Lam^0}({\bx}^0,{\by}^0)$ is obtained via
a functional ${q}^{\Lam^0}(\obom^0)=q^{\Lam^0}_{\vphi ,\Gam} 
(\obom^0)$ referred to as an infinite-volume RDMF: 
$$F^{\Lam^0}({\bx}^0,
{\by}^0)=\diy\int\limits_{\oW_{\bx^0,\by^0}}
{\bbP}^{\beta}_{{\bx}^0,{\by}^0}({\rd}\obom^0){q}^{\Lam^0}(\obom^0).
\eqno (3.1.5)$$
Further, the infinite-volume RDMF 
for a state $\vphi$ under consideration, should   
admit a particular representation. Namely, there exists a probability measure 
$\mu_\Gam =\mu_{\vphi,\Gam}$ on $(W_\Gam ,\fW_\Gam)$ such that 
$\forall$ finite set $\Lam'\subset\Gam$ with $\Lam^0\subset\Lam'$, 
$${q}^{\Lam^0}(\obom^0)=
\int_{W_{\Gam\setminus\Lam'}}
{\rd}\mu_{\Gam}^{\Gam\setminus\Lam'}
(\bom_{\Gam\setminus\Lam'})
{q}^{\Lam^0|\Gam\setminus\Lam'}(\obom^0|
\bom_{\Gam\setminus\Lam'})\eqno (3.1.6)$$
where
$${q}^{\Lam^0|\Gam\setminus\Lam'}(\obom^0|
\bom_{\Gam\setminus\Lam'})=
\frac{\Xi_{\Lam'\setminus\Lam^0}
(\obom^0,\bom_{\Gam\setminus\Lam'})}{\Xi_{\Lam'}
(\bom_{\Gam\setminus\Lam'})}.\eqno (3.1.7)$$
Here $\mu_{\Gam}^{\Gam\setminus\Lam'}$ stands for the 
restriction of measure $\mu_{\Gam}$ to the sigma-algebra 
${\fW}_{\Gam\setminus\Lam'}$.

Moreover, the expressions $\Xi_{\Lam'\setminus\Lam^0}
(\obom^0,\bom_{\Gam\setminus\Lam'})$ and $\Xi_{\Lam'}
(\bom_{\Gam\setminus\Lam'})$ represent, as before, partition 
functions in $\Lam'\setminus\Lam^0$ and $\Lam'$, with
the corresponding boundary conditions:
$$\begin{array}{r}\Xi_{\Lam'\setminus\Lam^0}
(\obom^0,\bom_{\Gam\setminus\Lam'})=
\diy\int\limits_{W_{\Lam'\setminus\Lam^0}}{\rd}
\nu_{\Lam'\setminus\Lam^0}(
\bom_{\Lam'\setminus\Lam^0})\qquad\qquad{}\\
\times\exp\big[-h^{\Lam'}
(\obom^0\vee\bom_{\Lam'\setminus\Lam^0}|
\bom_{\Gam\setminus\Lam'})\big]\end{array}\eqno (3.1.8)$$
and
$$\begin{array}{r}
\Xi_{\Lam'}(\bom_{\Gam\setminus\Lam'})=
\diy\int\limits_{W_{\Lam'}}{\rd}
\nu_{\Lam'}(\bom_{\Lam'})
\exp\big[-h^{\Lam'}
(\bom_{\Lam'}|\bom_{\Gam\setminus\Lam'})\big].
\end{array}\eqno (3.1.9)$$
The functional $h^{\Lam'}$
is defined by formulas similar to (2.3.8), (2.3.9):
$$h^{\Lam'}(\bom_{\Lam'}|\bom_{\Gam\setminus\Lam'})=
h^{\Lam'}(\bom_{\Lam'})+ h
(\bom_{\Lam'}|\,|\bom_{\Gam\setminus\Lam'})\eqno (3.1.10)$$
and
$$\begin{array}{r}
h^{\Lam'}(\obom^0\vee\bom_{\Lam'\setminus\Lam^0}|\bom_{\Gam\setminus\Lam'})=
h^{\Lam'}(\obom^0\vee\bom_{\Lam'\setminus\Lam^0})\qquad\qquad{}\\
+h
(\obom^0\vee\bom_{\Lam'\setminus\Lam^0}|\,|\bom_{\Gam\setminus\Lam'})
\end{array}\eqno (3.1.11)$$
where
$$h(\bom_{\Lam'}|\,|\bom_{\Gam\setminus\Lam'})=\sum\limits_{(j,j')\in
\Lam'\times (\Gam\setminus\Lam')}h^{j,j'}(\om_j,\om_{j'})
\eqno (3.1.12)$$
and similarly for $h^{\Lam'}(\obom^0\vee
\bom_{\Lam'\setminus\Lam^0}|\,|\bom_{\Gam\setminus\Lam'})$.
In turn, the terms $h^{\Lam'}$ 
and $h^{j,j'}$ are as in (2.3.8), (2.3.9). It is assumed that the series in
(3.1.12) is convergent for $\mu^{\Gam\setminus\Lam'}$-almost all
$\bom_{\Gam\setminus\Lam'})\in W_{\Gam\setminus\Lam'}$.

The functionals $h^{\Lam'}(\om_{\Lam'})$, 
$h^{\Lam'}(\bom_{\Lam'}|\bom_{\Gam\setminus\Lam'})$,
$h^{\Lam'}(\obom^0\vee\bom_{\Lam'\setminus\Lam^0})$,\\ $h^{\Lam'}(\obom^0\vee\bom_{\Lam'\setminus\Lam^0}|\bom_{\Gam\setminus\Lam'})$
and $h^{\Lam'}(\obom^0\vee
\bom_{\Lam'\setminus\Lam^0}|\,|\bom_{\Gam\setminus\Lam'})$ have the same 
meaning in terms of `energies' of loop/path configurations as before.

The measure $\mu_\Gam$ figuring in Eqns (3.1.6) and (3.1.8) has to satisfy the
infinite-volume DLR equations similar to (2.3.5). 
Namely, consider $p_{\Gam}^{\Lam^0}(\bom^0)$, the probability 
density function, relative to ${\rd}\nu^{\Lam^0}(\bom^0)$, 
for the measure $\mu_{\Gam}^{\Lam^0}$, the restriction to 
the sigma-algebra ${\fW}(\Lam^0)$ of measure $\mu_{\Gam}$:
$$p_{\Gam}^{\Lam^0}(\bom^0)=
\frac{{\rd}\mu_{\Gam}^{\Lam^0}(
\bom^0)}{{\rd}\nu^{\Lam^0}(\bom^0)}\,,\;\;
\bom^0\in W_{\Lam^0}\eqno (3.1.13)$$
The equations for $\mu_\Gam$ are that  
$\forall$ finite sets $\Lam^0$ and $\Lam'$ where
$\Lam^0\subset\Lam'\subset\Gam$, 
$$p_{\Gam}^{\Lam^0}(\bom^0)=
\int_{W_{\Gam\setminus\Lam'}}
{\rd}\mu_{\Gam}^{\Gam\setminus\Lam'}(\bom_{\Gam\setminus\Lam'})
p_\Gam^{\Lam^0|\Gam\setminus\Lam'}(\bom^0|
\bom_{\Gam\setminus\Lam'}).\eqno (3.1.14)$$
Here $p_\Gam^{\Lam^0|\Gam\setminus\Lam'}(\bom^0|
\bom_{\Gam\setminus\Lam'})$ is the conditional probability
density for $\bom^0$, conditioned by boundary condition
$\bom_{\Gam\setminus\Lam'}\in W_{\Gam\setminus\Lam'}$:    
$$p_\Gam^{\Lam^0|\Gam\setminus\Lam'}
(\bom^0|
\bom_{\Gam\setminus\Lam'})
=\frac{\Xi_{\Lam'\setminus\Lam^0}
(\bom^0,\bom_{\Gam\setminus\Lam^0})}
{\Xi_{\Lam'}(\bom_{\Gam\setminus\Lam'})}\,.\eqno (3.1.15)$$
Here, as in (3.1.8)
$$\begin{array}{r}\Xi_{\Lam'\setminus\Lam^0}
(\bom^0,\bom_{\Gam\setminus\Lam^0})=
\diy\int\limits_{W_{\Lam'\setminus\Lam^0}}{\rd}
\nu_{\Lam'\setminus\Lam^0}(
\bom_{\Lam'\setminus\Lam^0})\qquad{}\\
\times\exp\big[-h^{\Lam'}
(\bom^0\vee\bom_{\Lam'\setminus\Lam^0}|
\bom_{\Gam\setminus\Lam'})\big]\end{array}\eqno (3.1.16)$$
and, as in (3.1.9), 
$$\begin{array}{r}
\Xi_{\Lam'}(\bom_{\Gam\setminus\Lam'})=
\diy\int\limits_{W_{\Lam'}}{\rd}
\nu_{\Lam'}(\bom_{\Lam'})
\exp\big[-h^{\Lam'}
(\bom_{\Lam'}|\bom_{\Gam\setminus\Lam'})\big],
\end{array}\eqno (3.1.17)$$
with the functional $h^{\Lam'}$
defined similarly to Eqns (3.1.10)--(3.1.12).
\vskip 5 truemm

{\bf Remark 3.1.} We do not claim (at least in this paper and its sequel
\cite{KS2}) that
the properties (3.1.4)--(3.1.17) imply that the operator $R^{\Lam^0}$
is an RDM (positive definiteness of $R^{\Lam^0}$ remains an 
open question). However, when one can assert (on grounds of some additional
information) that a given family of operator $\{R^{\Lam^0}\}$ obeying
(3.1.4)--(3.1.17) consists of positive-definite operators then we can speak of
a state $\vphi\in\fG$. (Viz., this is the case of limiting Gibbs states
discussed in Theorems 1.1 and 1.2.) As it stands, our results stated in
Section 3.2 hold true for any family of operators $R^{\Lam^0}$ for which
Eqns (3.1.4)--(3.1.17) are fulfilled. E.g., we can claim the assertion
of Theorem 1.2 for any linear normalized functional on $\fB$ defined by a family
$\{R^{\Lam^0}\}$ satisfying (3.1.4)--(3.1.17).  
\vskip 5 truemm

Some elements of the above construction have been used in the literature;
see, e.g., \cite{KoP} and references therein.

We will refer to Eqns (3.1.6)--(3.1.12) as an FK-DLR representation (of 
the infinite-volume RDMF ${q}^{\Lam^0}(\omega^0 )$) by a given probability 
measure $\mu_\Gam$, assuming that $\mu_\Gam$ satisfies the
infinite-volume DLR equations (3.1.14). It is important  
to stress that, unlike the case of a finite $\Lam\subset\Gam$,
the solution to the infinite-volume FK-DLR equations (3.1.14)  
over the whole graph $\Gam$
may be, in general, non-unique. However, the family of functionals
${q}^{\Lam^0}$, where $\Lam^0$ runs through the finite 
subsets of $\Gam$ 
is determined uniquely provided that a measure $\mu_\Gam$
is given, satisfies (3.1.6), (3.1.7). In accordance with
the above scheme, this gives rise to the family 
of RDMKs $F^{\Lam}$ and -- ultimately -- RDMs $R^{\Lam}$, 
for finite sets $\Lam\subset\Gam$, obeying the above compatibility
property (1.4.9). The corresponding state 
(emerging from the probability  measure 
$\mu_{\Gam}$) is denoted by $\vphi_{\Gam} (=
\vphi (\mu_{\Gam}))$; when possible, the subscript $\Gam$
will be omitted. Given $\beta\in (0,\infty )$, the class of the 
measures $\mu_{\Gam}=\mu_{\beta ,\Gam}$
satisfying Eqns (3.1.13) is denoted by ${\fG}(\beta )$, 
as well as the class of related states $\vphi_{\Gam}$.    

In Theorem 3.1 below we establish that the 
class $\fG (\beta )$ is non-empty $\forall$ given $\beta\in (0,\infty )$.

As was said earlier, the infinite-volume invariance property 
under study is expressed by Eqn (3.1.1):
$\forall$ ${\ttg}\in{\ttG}$, 
finite set $\Lam^0\subset\Gam$, $\bx^0,\by^0\in M^{\Lam^0}$ 
and $\obom^0=\{\oomega_j,
\;j\in\Lam^0\}\in\oW_{\bx^0,\by^0}$, the value 
${q}_{\beta}^{\Lam^0}({\ttg}\obom^0)=
{q}_{\beta}^{\Lam^0}(\obom^0)$. Here ${\ttg}\obom^0$ is as in (2.3.16).
A similar property for the density $p^{\Lam^0}
(\obom^0)$ has the form: $\forall$
finite set $\Lam^0\subset\Gam$ and loop configuration
$\obom^0=\{\oomega_j,\;j\in\Lam^0\}\in \oW_{\Lam^0}$, 
$${p}^{\Lam^0}(\obom^0)={p}^{\Lam^0}(
{\ttg}\obom^0),\;\;{\ttg}\in{\ttG}.\eqno (3.1.18)$$

The invariance properties in Eqns (3.1.1) and (3.1.18) imply that, $\forall$
finite set $\Lam^0\subset\Gam$, the 
infinite-volume RDMs $R_\Gam^{\Lam^0}$ have the property similar to (1.4.9):
$$R_\Gam^{\Lam^0}U_{\Lam^0}({\ttg})=
U_{\Lam^0}({\ttg})R_\Gam^{\Lam^0}\eqno{(3.1.19)}$$
which, in terms of the corresponding state $\vphi$, means 
(1.4.10).
\vskip .5 truecm

{\bf 3.2. Properties of class $\fG$.} The results of this 
paper about state class $\fG$  are summarised in 
Theorems 3.1--3.4 below.
In Theorems 3.1--3.3 we assume the above conditions (1.1.1), (1.1.2), 
(1.3.2)--(1.3.5).
\vskip .5 truecm

{\bf Theorem 3.1.} {\sl 
For all $\beta\in (0,\infty )$, the sequence of 
Gibbs states $\vphi_{\Lam (n)}$ contains a subsequence 
$\vphi_{\Lam(n_k)}$
such that $\forall$ finite $\Lam^0\subset\Gam$ and $A_0\in
\fB_{\Lam^0}$, we have:
$$\lim_{k\to\infty}\vphi_{\Lam(n_k)}(A_0)=\vphi (A_0)$$
where state $\vphi\in\fG (\beta )$. Consequently, class $\fG (\beta )$ 
is non-empty.$\quad\lhd$} 
\vskip .5 truecm

{\bf Theorem 3.2.} {\sl For all
$\beta\in (0,\infty )$ and finite $\Lam^0\subset\Gam$, any Gibbs 
state $\vphi\in{\fG}(\beta )$ satisfies properties {\rm{(3.1.1)}} and 
{\rm{(3.1.18)--(3.1.19)}}.$\quad\lhd$}
\vskip .5 truecm

The invariance property can be formally extended to ground states.
We call a state $\ovphi$ (of C$^*$-algebra $\fB$) a ground state 
if there exists a sequence of states
$\vphi_n\in{\fG}(\beta_n)$ with $\beta_n\to\infty$ such that 
$\ovphi ={\rm{w}}^*-\lim\limits_{n\to\infty}\vphi_n$. 
\vskip .5 truecm

{\bf Corollary 3.3.} {\sl Any ground state $\ovphi$
(i.e., a w$^*$-limiting point of states $\vphi_n\in\fG(\beta_n)$)
obeys {\rm{(3.1.1)}} and {\rm{(3.1.18)--(3.1.19)}}.$\quad\lhd$}
\vskip .5 truecm

Of course, Corollary 3.3 does not prove existence of ground states for the 
model under consideration.

In a future paper we will remove the smoothness condition 
upon the potential
function $V$ (see (1.3.2)), by following the methodology from 
\cite{R1}. 
In this paper we note that, like the classical case (cf. \cite{ISV} 
and the bibliography therein), if the condition of smoothness is 
violated, the symmetry property may be destroyed. See Theorem 3.4 below. 
\vskip .5 truecm

{\bf Theorem 3.4.} {\sl Take $\Gam ={\bbZ}^2$, the regular square
lattice, with distance ${\tt d}(j,j')=\max\,\big[|j_1-j'_1|,
|j_2-j'_2|\big]$. 
Take $M=S^1={\ttG}$ where $S^1={\bbR}/{\bbZ}$ is 
a unit circle, with a standard metric $\rho (x,x')=\min\,\big[|x-x'|,
1-|x-x'|\big]$ and 
the group operation of addition mod $1$. Assume 
that the two-body potentials $J({\ttd}(j,j'))$ and $V(x,x')$, $j,j'\in\bbZ^2$,
$x,x'\in S^1$, 
are of the form 
$$\begin{array}{c}J({\ttd}(j,j'))=\begin{cases}1,&|j-j'|=1,\\ 0,&|j-j'|\neq 1,\end{cases}\\
V(x,x')=\begin{cases}-\cos\,2\pi(x-x'),&
\rho (x,x')\leq\theta,\\
+\infty,\rho (x,x')>\theta,
\end{cases}\end{array}\eqno (3.2.1)$$
with a usul agreement $0\cdot (+\infty )=0$, where $\theta\in (0,1/4)$ is
a constant.  
In this case, Hamiltonian $H_\Lam$ is
equipped with Dirichlet's boundary conditions on $D\subset M^{\Lam}$
where 
$$D=\Big\{\bx_\Lam\in M^{\Lam}:\;\hbox{ $|x(j)-x(j')|\geq\theta$
for some $j,j'\in\Lam$ with ${\ttd}(j,j')=1$}\Big\}.$$
Then, $\forall$ $\beta\in (0,\infty )$, there exists an FK-DLR measure 
${\wt\mu}={\wt\mu}_\beta\in\fG$ which is not $S^1$-invariant. Consequently, the 
corresponding FK-DLR state
${\wt\vphi}={\wt\vphi}_{\wt\mu}\in{\fG}(\beta )$ is not 
$S^1$-invariant.$\quad\lhd$}
\vskip .5 truecm

Similarly to a ground state $\ovphi$, we can define $\omu$, a ground-state
FK-measure. Namely, take an FK-DLR measure $\mu_n\in\fG(\beta_n)$ 
and consider its image ${\wh\mu}_n$ under projection 
$\bom_\Gam\in W^\beta_\Gam\mapsto \bx_\Gam
\in M^\Gam$ where $\bx_\Gam =\bx_\Gam (\bom_\Gam )$ is the 
collection of initial points for loop configuration $\bom_\Gam$. 
Suppose that ${\wh\mu}$
is a limiting point for sequence ${\wh\mu}_n$ as $\beta_n\to\infty$. Then 
we say that ${\wh\mu}$ is a ground-state FK-measure. Furthermore, such
a measure is called ${\ttG}$-invariant if $\forall$ finite set 
$\Lam^0\subset\Gam$, ${\ttg}\in{\ttG}$ and a (bounded)
function $\phi:M^{\Lam^0}\to\bbR$, the integral 
${\wh\mu}(U_{\Lam^0}({\ttg})\phi)={\wh\mu}(\phi)$.
\vskip .5 truecm

{\bf Corollary 3.5.} {\sl Suppose that the family of non-$S^1$-invariant 
FK-DLR measures ${\wt\mu}_\beta$ specified in Theorem {\rm{3.4}}
has a limiting point ${\wt\psi}$ as $\beta\to\infty$. Then ${\wt\psi}$
is a non-$S^1$-invariant ground-state FK-measure.$\quad\lhd$}  
\vfill\eject

\centerline{\bf 4. Proof of the main results}
\vskip 1 truecm

In this section we deliver proofs of the stated results.
\vskip .5 truecm

{\bf 4.1. Proof of Theorem} 3.1. Given finite sets $\Lam^0$ and $\Lam$, 
$\Lam^0\subset\Lam\subset\Gam$,
the RDMK $F^{\Lam^0}_\Lam$ (see (2.2.1)) is a continuous function on 
$M^{\Lam^0}\times M^{\Lam^0}$. The first observation is that 
the sequence $F^{\Lam^0}_{\Lam (n)}$ is compact in 
$C^0(M^{\Lam^0}\times M^{\Lam^0})$ by the Ascoli--Arzela theorem,
as these functions are uniformly bounded and continuous. The latter 
property is based upon conditions (1.3.2)--(1.3.4) and the assumption 
that $M$ is compact.

More precisely, to show uniform boundedness, note that $\forall$
finite $\Lam'\subset\Gam$ with $\Lam'\supset\Lam^0$, the conditional
RDMF ${q}^{\Lam^0|\Gam\setminus\Lam'}$ (cf. Eqn (3.1.7)) satisfies
$${q}^{\Lam^0|\Gam\setminus\Lam'}(\obom^0|
\bom_{\Gam\setminus\Lam'})\leq (e^{2\beta\oJ{\ov V}})^{\sharp\,\Lam^0}
\eqno (4.1.1)$$ 
for all path configurations $\obom^0\in \oW_{\Lam^0}$,
$\bom_{\Gam\setminus\Lam'}\in W_{\Gam\setminus\Lam'}$. Here, the constant $\oJ$ is given by
$$\oJ=\sup\left[\sum\limits_{j'\in\Gam}J({\tt d}(j,j')):\;j\in\Gam \right]\eqno (4.1.2)$$
where $J$ is the function from (1.3.3). (In fact, ${\ov J}$ coincides with the quantity
${\ov J}(1)$ in (1.3.1).) A similar bound holds if we
replace $\Gam$ with ball $\Lam (n)\subset\Gam\setminus\Lam'$.
After integration, this yields the estimate
$$F^{\Lam^0}_{\Lam (n)}(\bx^0,\by^0)\leq 
(e^{2\beta\oJ{\ov V}}\,{\ov p}^{\beta})^{\sharp\,\Lam^0}\eqno (4.1.3)$$
where 
$${\ov p}^{\beta}=\sup\big[p^\beta (x,y), |\nabla_xp^\beta (x,y)|,
|\nabla_yp^\beta (x,y)|:\;x,y\in M\big].\eqno (4.1.4)$$
Here $p^\beta (x,y)$ stands for the transition probability density
from $x$ to $y$ for the Brownian motion (with the generator 
$-\Delta/2$ where $\Delta$ is the Laplace operator on $M$)
in time $\beta$:
$$p^\beta (x,y)=\frac{1}{(2\pi\beta)^{d/2}}\sum_{\un =(n_1,\ldots ,n_d)\in\bbZ^d} 
\exp\,\Big(-|x-y+\un |^2/2\beta\Big)\,,\eqno (4.1.5)$$ 

The argument for uniform continuity (or equi-continuity) of RDMKs
$F^{\Lam^0}_{\Lam (n)}(\bx^0,\by^0)$ is more technical. We want to 
check that the gradients\\ $\nabla_{\bx^0}F^{\Lam^0}_{\Lam (n)}(
\bx^0,\by^0)$ and $\nabla_{\by^0}F^{\Lam^0}_{\Lam (n)}(\bx^0,\by^0)$
are uniformly bounded. There are two contributions into the gradient:
one comes from varying the measure $\bbP_{\bx^0,\by^0}({\rd}\obom^0)$, the other from
varying the functional\\ $\exp\,[-h^{\Lam^0}
(\obom^0\,|\,\bom_{\Lam (n)\setminus\Lam^0})]$. The first
contribution can be
uniformly bounded in terms of the constant ${\ov p}^{\beta}$.

The second contribution can be analysed by deforming a path $\oom_j\in\obom^0$,
$j\in\Lam^0$: one of the end-points $x(j)$ or $y(j)$ can be moved, say,  
along a geodesic on $M$ (i.e., a straight line). The 
points on the path are then moved, at a scaled distance, via a parallel 
transfer (the affine connection
on the torus is trivial). This contribution
is related to the differentiation of the exponent and
controlled due to the bound (1.3.2) on the derivatives of the potential.  

The second contribution yields an expression of the form
$$\begin{array}{r}\diy\int_{\oW_{\bx^0,\by^0}}
\bbP_{\bx^0,\by^0}({\rd}\obom^0)\sum_{j\in\Lam^0}
{\wth}_j(\om (j),\bom_{\Lam (n)\setminus\Lam^0})\qquad{}\\    
\diy\times\exp\,\left[-h^{\Lam^0}
(\obom^0\,|\,\bom_{\Lam (n)\setminus\Lam^0})\right]
\end{array}\eqno (4.1.6)$$
where functional ${\wth}_j(\obom^0,\bom_{\Lam (n)\setminus\Lam^0})$
is uniformly bounded. Combining this with the argument used to estimate
the RDMK $F^{\Lam^0}_{\Lam (n)}(\bx^0,\by^0)$ allows one to    
bound the gradients $\nabla_{\bx^0}
F^{\Lam^0}_{\Lam (n)}(\bx^0,\by^0)$ and 
$\nabla_{\by^0}F^{\Lam^0}_{\Lam (n)}(\bx^0,\by^0)$ as well.

More precisely, given a site $j\in \Lambda^0$,
we need to differentiate in $\ux(j)$ or $\uy(j)$ the expression
$$
\begin{array}{l}\diy
\sum\limits_{\bn}\int_{\oW_{\bx^0,\by^0+\bn}}{\bbP}_{\bx^0,\by^0+\bn}(d\obom)
\exp[-h^{\Lam^0}(\obom|
\bom_{\Lam(n)\setminus \Lam^0}]\\
\diy\quad =\sum\limits_{\bn} \exp\;\Big(-|\bx^0-\by^0-\bn|^2/2\beta\Big)\\
\diy\qquad\times\int_{W_{\bx^0,\bx^0}}{\bbP}_{\bx^0,\bx^0}(d\obom^0)
\exp\Big[-h^{\Lam^0}(\obom^0+\zeta_{\bn}|\bom_{\Lam(n)\setminus\Lam^0})
\Big].\end{array}$$ 
Here we sum over vectors $\bn =(\un (j),j\in\Lam^0)\in (\bbZ^d)^{\Lam^0}$. 
Furthermore, $\zeta_{\bn}$ is a linear map:
$$\zeta_{\bn}(\tau)=\diy\frac{\tau}{\beta}(\by^0+\bn-\bx^0),\;\;
0\leq\tau\leq\beta .$$
Finally, the measures ${\bbP}_{\bx^0,\by^0+\bn}$ and ${\bbP}_{\bx^0,\bx^0}$ 
refer to the standard Brownian motion on $\bbR^d$.
 
Suppose we differentiate in $y(j)\in\bbR^d$. Differentiating 
the exponent yields a convergent series. Next, 
$$\begin{array}{r}
\diy\nabla_{y(j)}\exp\Big[-h^{\Lam^0}(\obom^0
+\zeta_{\bn}|\bom_{\Lam(n)\setminus\Lam^0})\Big]\qquad\qquad\qquad\qquad{}\\
\diy =\left(\sum\limits_{l\in\Lam (n)\setminus\Lam^0}
\int_0^\beta{\rm d}\tau \frac{\tau}{\beta}(\nabla_{y(j)}V)
(\om (j,\tau ),\om (l,\tau))\right)\quad{}\\
\diy\times\exp\Big[-h^{\Lam^0}(\obom^0
+\zeta_{\bn}|\bom_{\Lam(n)\setminus\Lam^0})\Big]
\end{array}$$
is bounded due to (1.3.2). (The expression in the big brackets gives 
the term ${\wth}_j(\om (j),\bom_{\Lam (n)\setminus\Lam^0})$ figuring in 
Eqn (4.1.6).) This yields the desired result.

Differentiation in $x(j)$ can be done in a similar manner, by exchanging 
$\bx^0$ and $\by^0$ in the above series.

Now let an RDMK $F^{\Lam^0}$ be a limiting point for
$F^{\Lam^0}_{\Lam (n_k)}$ in $C^0(M^{\Lam^0}\times M^{\Lam^0})$.
Then we have that 
$$\diy\lim_{k\to\infty}
\int\limits_{M^\Lam\times M^\Lam}{\rd}\bx^0{\rd}\by^0
\left[F^{\Lam^0}_{\Lam (n_k)}(\bx^0,\by^0)-F^{\Lam^0}(\bx^0,\by^0)
\right]^2=0.$$
In other words, the RDM $R_{\Lam (n_k)}^{\Lam^0}$ in $\cH^{\Lam^0}$
converges to the infinite-volume RDM $R^{\Lam^0}$ in the 
Hilbert-Schmidt norm: $\left\|R_{\Lam (n_k)}^{\Lam^0}
-R^{\Lam^0}\right\|_{\rm{HS}}\to 0$. According to Lemma 1 from \cite{Su1},
the convergence takes place in the trace-norm as well:
$\left\|R_{\Lam (n_k)}^{\Lam^0}
-R^{\Lam^0}\right\|_{\rtr}\to 0$. We obtain that the sequence of states
$\vphi_{\Lam (n)}$ is w$^*$-compact. 

In parallel, an argument can be developed that the measures 
$\mu_\Lam$ form a compact family as $\Lam\nearrow\Gam$.
More precisely, we would like to show that $\forall$ 
finite $\Lam^0\subset\Gam$, the family of measures $\mu^{\Lam^0}_\Lam$
is compact. To see this, it suffices to check that,
for a fixed $\Lam^0$,  the sequence
$\{\mu^{\Lam (L)}_{\Lam (n)},\;n=L+1,L+2,\ldots\}$
is tight and apply the Prokhorov theorem.

To check tightness, we use the two facts. (i) The reference measure
${\rd}\nu_{\Lam^0}$ on $W^{\Lam^0}$ (see Eqn (2.3.1)) is supported by 
loop configurations
with the standard continuity modulus $\sqrt{2\epsilon\ln\,(1/\epsilon)}$.
(ii) The probability density 
$\diy p_\Lam^{\Lam^0}(\om^0)=\frac{{\rd}\mu^{\Lam (L)}_{\Lam (n)}(\Om (\Lam (L)))}{
{\rd}\nu_{\Lam (n)}(\Om (\Lam (L)))}$ (cf. Eqn (2.3.3))
is bounded from above by a constant $\exp \,\big\{\beta
\big[\sharp\Lam (L)\big]{\ov V}{\ov J}^*\big\}$
(and from below by $\exp \,\big\{-\beta
\big[\sharp\Lam (L)\big]{\ov V}J^*\big\}$). See (1.3.2) and (1.3.4).

The constructed family of limit-point measures $\mu^{\Lam^0}$
has the compatibility property and therefore satisfies the assumptions
of the Kolmogorov theorem. The result is that there exists a unique
probability measure $\mu$ on $W^\Gam$ such that the restriction of $\mu$ on
the sigma-algebra of subsets localized in $\Lam^0$ coincides with
$\mu^{\Lam^0}$. 

The fact that $\mu$ is FK-DLR follows from the above construction.
Hence, each limit point $\vphi$
falls in class $\fG (\beta )$. $\quad\Box$
\vskip .1truecm
 
{\bf Remark 4.1.} Anticipating a forthcoming result for a general compact 
manifold $M$, we propose to discuss a version of the above argument
for an example where $M$ is a two-dimensional Klein bottle
with a flat Riemannian metric. A convenient representation is through 
the universal simply connected cover which in this case is the Euclidean 
plane $\bbR^2$ with the standard metric. For the fundamental polygon we take
a square $[-1/2,1/2]\times [-1/2,1/2]$ where the following 
pairs of points $({\rx}_1;{\rx}_2)$ are glued: 
$$({\rx}_1;-1/2)\;\hbox{ and }\;(-{\rx}_1;1/2),\;\hbox{ where }\;-1/2<{\rx}_1<1/2$$ 
and
$$(-1/2;{\rx}_2)\;\hbox{ and }\;(1/2;{\rx}_2),\;\hbox{ where }\;-1/2<{\rx}_2<1/2.$$
The cover map $T:\,\bbR^2\to M$ is as follows:
$$\begin{array}{c}T ({\rx}_1,{\rx}_2)\mapsto ((-1)^{n_2}({\rx}_1-n_1);{\rx}_2-n_2)\\
\hbox{ whenever }\;-1/2\leq{\rx}_i-n_1<1/2,\;-1/2\leq{\rx}_2-n_2<1/2.
\end{array}\eqno (4.1.7)$$
Here and below, $n_1,n_2\in\bbZ$ are integers.
(In this example, ${\ttG}$ may be a circle ${\tt S}^1$ (a 1D
torus) realized as an interval $[-1/2,1/2]$ with points $-1/2$ and $1/2$ glued
together. The action  is: ${\ttg}x=({\rx}_1+{\ttg},{\rx}_2)$, for $x=({\rx}_1;{\rx}_2
)$, with addition in $[-1/2,1/2]$. (Other choices of ${\tt G}$ (a 1D torus 
of length 2 or a 2D torus) are analysed in a similar fashion.)

In this example, the integral 
$$\int_{\oW_{\bx^0,\by^0}}\bbP_{\bx^0,\by^0}({\rd}\obom^0)\exp\,\big[-
h^{\Lam^0}(\obom^0|
\bom_{\Lam (n)\setminus\Lam^0})\big]\eqno (4.1.8)$$
contributing to $F^{\Lam^0}_\Lam (\bx^0,\by^0)$ can again be 
differentiated explicitly. For definiteness, take $\Lam^0$
to be a one-point set, say $\{o\}$, with particle configurations $\bx^0$
and $\by^0$ reduced to single points in $M$ (or rather in the 
fundamental polygon): 
$$\bx^0=x=({\rx}_1;{\rx}_2),\; \by^0=y=({\ry}_1;{\ry}_2),\;-1/2\leq{\rx}_i,
{\ry}_i\leq 1/2,\;i=1,2.$$ 
Then the above integral takes the form
$$\begin{array}{l}\diy\frac{1}{2\pi\beta}\sum_{n_1,n_2\in\bbZ}
\exp\,\left[-((-1)^{n_2}{\rx}_1-{\ry}_1-n_1)^2/(2\beta)\right]\\
\qquad\qquad\times\exp\,\left[-({\rx}_2-{\ry}_2-n_2)^2/(2\beta)\right]
\diy\int\limits_{W^{\beta ,\bbR^2}_{x,x}}\bbP_{x,x}^{\beta ,\bbR^2}
({\rd}\om^0_1\times{\rd}\om^0_2)\\
\qquad\times\exp\left[-h^{\Lambda^0}
\left(T\left[\left(\om^0_1+\delta_1^{(n_1)},\om^0_2
+\delta_2^{(n_2)}\right)\right]\Big|\bom_{\Lam (n)\setminus\Lambda^0}
\right)\right].\end{array}\eqno (4.1.9)$$  
Here $W^{\beta ,\bbR^2}_{x,x}$ and $\bbP_{x,x}^{\beta ,\bbR^2}$ stand, 
respectively, for the space of continuous loops (closed trajectories 
beginning and ending at $x$) and the Brownian bridge distribution in 
the plane $\bbR^2$ of the time-length $\beta$. Next, $\delta_i^{(n_i)}$ 
are functions $[0,\beta ]\to\bbR$ providing deformations of plane loops
$\om^0_i$ from $W^{\beta ,\bbR^2}_{{\rx}_i,{\rx}_i}$ into plane paths
$\om^0_1+\delta_1^{(n_1)}$ from 
$\oW^{\beta ,\bbR^2}_{{\rx}_i,{\ry}_i+n_i}$:
$$\delta_1^{(n_1)}(\tau )=\frac{\tau}{\beta}\left({\ry}_1-{\rx}_1+n_1\right),\;\;
\delta_1^{(n_2)}(\tau )=\frac{\tau}{\beta}\left({\ry}_2-{\rx}_2+n_2\right).
\eqno (4.1.10)$$

A similar formula holds after replacing
$W^{\beta ,\bbR^2}_{x,x}$ and $\bbP_{x,x}^{\beta ,\bbR^2}$ with
$W^{\beta ,\bbR^2}_{y,y}$ and $\bbP_{y,y}^{\beta ,\bbR^2}$,
{\it mutatis mutandis}. Differentiations $\nabla_x$ and $\nabla_y$
then become straightforward (although rather tedious), confirming the above 
claim about uniform continuity of RDMKs  
$F^{\Lam^0}_\Lam $ (cf. Eqn (4.1.6)).
\vskip .5 truecm

Theorem 1.1 can be deduced from Theorem 3.1 since every limiting
Gibbs state $\vphi\in\fG^0$ lies in $\fG$ (in other words, $\fG^0\subseteq\fG$), 
as follows from the above argument and the observations made in Sections 2.3 and 3.1.
\vskip .5 truecm

{\bf 4.2. Proof of Theorem 3.2.} We follow  
the approach initiated in \cite{FP}. The proof of Theorem 3.2 is
based on the following bound for infinite-volume RDMFs: $\forall$ 
finite $\Lam^0\subset\Gam$, $\obom^0\in\oW_\Lam$ and 
${\ttg}\in{\ttG}$, 
$${q}^{\Lam^0}({\ttg}\obom^0)+{q}^{\Lam^0}({\ttg}^{-1}\obom^0)
\geq 2{q}^{\Lam^0}(\obom^0).\eqno (4.2.1)$$
Lemma 1 from \cite{FP} implies that (3.1.1) follows from (4.2.1).

In turn, the bound (4.2.1) follows from a similar inequality for the
conditional RDMFs: $\forall$ finite $\Lam^0\subset\Gam$, 
$\obom^0\in\oW_\Lam$, ${\ttg}\in{\ttG}$ and $a \in (1,\infty )$,
for any $n$ large enough and $\bom_{\Gam\setminus\Lam (n)}$,
$$\begin{array}{r}a{q}^{\Lam^0|\Gam\setminus\Lam (n)}({\ttg}\obom^0|
\bom_{\Gam\setminus\Lam (n)})
+a{q}^{\Lam^0|\Gam\setminus\Lam (n)}({\ttg}^{-1}\obom^0|
\bom_{\Gam\setminus\Lam (n)})\qquad{}\\
\geq 2{q}^{\Lam^0|\Gam\setminus\Lam (n)}(\obom^0|
\bom_{\Gam\setminus\Lam (n)}).\end{array}\eqno (4.2.2)$$
In fact, to deduce (4.2.1) from (4.2.2), it is enough to 
integrate (4.2.2) in\\ ${\rd}\mu_\Gam^{\Gam\setminus\Lam (n)}
(\bom_{\Gam\setminus\Lam (n)})$ and let $a \searrow 1$. 

Now, Eqn (4.2.2) is deduced after performing a special construction
related to a family of `gauge' actions $\bttg_{\Lam (n)\setminus\Lam^0}$
on loop configurations $\bom_{\Lam (n)\setminus\Lam^0}$; see Eqns (4.2.4), (4.2.5) below.
A particular feature of action $\bttg_{\Lam (n)\setminus\Lam^0}$ is that 
it `decays' to ${\tt e}$, the unit element of ${\ttG}$ (which generates
a `trivial' identity action), when we move from $\Lam^0$ towards $\Gam
\setminus\Lam (n)$. 
Formally, (4.2.2) will follow from the inequality: $\forall$ finite 
$\Lam^0\subset\Gam$, 
$\obom^0\in\oW_\Lam$, ${\ttg}\in{\ttG}$ and $a \in (1,\infty )$,
for any $n$ large enough, $\bom_{\Lam (n)\setminus\Lam^0}$ and 
$\bom_{\Gam\setminus\Lam (n)}$, 
$$\begin{array}{l}\diy\frac{a}{2}\exp\,\Big[-h^{\Lam (n)}(({\ttg}\obom^0)\vee
(\bttg_{\Lam (n)\setminus\Lam^0}\bom_{\Lam (n)\setminus\Lam^0})|
\bom_{\Gam\setminus\Lam (n)})\Big]\\
\qquad\quad\diy +
\frac{a}{2}\exp\,\Big[-h^{\Lam (n)}(({\ttg}^{-1}\obom^0)\vee
(\bttg^{-1}_{\Lam (n)\setminus\Lam^0}\bom_{\Lam (n)\setminus\Lam^0})|
\bom_{\Gam\setminus\Lam (n)})\Big]\\
\qquad\qquad\geq\diy
\exp\,\Big[-h^{\Lam (n)}(\obom^0\vee\bom_{\Lam (n)\setminus\Lam^0}|
\bom_{\Gam\setminus\Lam (n)})\Big]\,.
\end{array}\eqno (4.2.3)$$
Indeed, (4.2.2) follows from (4.2.3) by integrating in ${\rd}\nu_{\Lam (n)
\setminus\Lam^0}(\bom_{\Lam (n)\setminus\Lam^0})$ and normalizing
by $\Xi_{\Lam (n)\setminus\Lam^0}(\bom_{\Gam\setminus\Lam (n)})$;
cf. Eqn (3.1.16) with $\Lam'=\Lam (n)$.
(Here one uses the fact that the Jacobian of the map $\bom_{\Lam (n)
\setminus\Lam^0}\mapsto\bttg_{\Lam (n)\setminus\Lam^0}
\bom_{\Lam (n)\setminus\Lam^0}$ equals $1$.)

Thus, our aim becomes to prove (4.2.3). The gauge family 
$\bttg_{\Lam (n)\setminus\Lam^0}$ is composed by individual
actions ${\ttg}^{(n)}_j\in{\ttG}$: 
$$\bttg_{\Lam (n)\setminus\Lam^0}=\{{\ttg}^{(n)}_j,\;j\in\Lam (n)
\setminus\Lam^0\}.\eqno (4.2.4)$$
Let us identify the element 
${\ttg}\in{\tt G}$ with a vector $\utheta =\theta A\in M$ and use the additive notation:
${\ttg}x:=x+\utheta$, $x\in M$. Then ${\ttg}^{(n)}_j\in {\ttG}$ corresponds to multiples of
thye vector $\utheta$. Namely, we fix 
a positive
integer value $\ovr$ such that $\Lam^0\subset\Lam_{\ovr}$ and identify
$${\ttg}^{(n)}_j\;\hbox{ with \;\;}\utheta {\upsilon (n,j)}\eqno (4.2.5)$$
where 
$$\upsilon (n,j)
=\begin{cases}1,&{\tt d}(o,j)\leq\ovr,\\
\vartheta\big({\tt d}(j,o)-\ovr,n-\ovr\big),&{\tt d}(o,j)>\ovr .
\end{cases}\eqno (4.2.6)$$
In turn, the function $\vartheta (a,b)$ satisfies
$$\vartheta (a,b)={\mathbf 1}(a\leq 0)+
\frac{{\mathbf 1}(0<a<b)}{Q(b)}\int_a^bz(u){\rd}u,\;\;
a,b\in\bbR ,\eqno (4.2.7)$$
with the same functions $Q(b)$ and $z(u)$ as proposed in \cite{FP}
$$\begin{array}{l}\diy Q(b)=\int_0^b z(u){\rd}u,\\
\diy\quad\hbox{where}
\;z(u)={\mathbf 1}(u\leq 2)+{\mathbf 1}(u>2)
\frac{1}{u\ln\,u},\;b>0.\end{array}\eqno (4.2.8)$$
Moreover, $\bttg_{\Lam (n)\setminus\Lam^0}^{-1}$ is the 
collection of the inverse elements:
$$\bttg_{\Lam (n)\setminus\Lam^0}^{-1}=
\left\{{{\ttg}_j^{(n)}}^{-1},\;j\in\Lam (n)\setminus\Lam^0\right\}.$$ 
We will use the formulas for ${\ttg}^{(n)}_j$ for $j\in\Lam (n)$, or even for
$j\in\Gam$,
as they agreed with the requirement that ${\ttg}^{(n)}_j\equiv{\ttg}$
when $j\in\Lam^0$ and ${\ttg}^{(n)}_j\equiv{\tt e}$ for $j\in\Gam
\setminus\Lam (n)$. Accordingly, we will use the notation $\bttg_{\Lam (n)}
=\{{\ttg}^{(n)}_j,j\in\Lam (n)\}$.

Next, we use the invariance property (1.3.5). The 
Taylor formula for function $V\in C^2$ yields for
$j,j'\in\Lam (n)$:
$$\begin{array}{l}
\Big|V\left({\ttg}^{(n)}_j\om_j,{\ttg}^{(n)}_{j'}\om_{j'}\right)\\
\qquad +V\left({{\ttg}^{(n)}_j}^{-1}\om_j,
{{\ttg}^{(n)}_{j'}}^{-1}\om_{j'}\right)
-2V(\om_j,\om_{j'})\Big|\\ \;\\
\qquad\qquad\leq  C\,
|\utheta|^2\left|\upsilon (n,j)-\upsilon (n,j')\right|^2
{\ov V}\,.\end{array}\eqno (4.2.9)$$ 
Here $C\in (0,\infty )$ is a constant ${\ov V}$ is taken from (1.3.2) and notations from
(4.2.5) are used. 

The bound (4.2.9) is crucial and exploits the structure of the group action. 
It uses the fact that the first-order terms
in the expansion in the left-hand side of (4.2.9) cancel each other
because of the presence of elements ${\ttg}^{(n)}_j$ and  
${\ttg}^{(n)}_{j'}$ and their inverses, ${{\ttg}^{(n)}_j}^{-1}$ and
${{\ttg}^{(n)}_{j'}}^{-1}$. (This idea goes back to \cite{P} and \cite{FP}.)

The term
$|\upsilon (n,j)-\upsilon (n,j')|^2$ can be specified as 
$$|\upsilon (n,j)-\upsilon (n,j')|^2=\begin{cases}0,
\;\hbox{ if }\;{\tt d}(j,o),
{\tt d}(j',o)\leq \ovr,\\
0,\;\hbox{ if }\;{\tt d}(j,o),{\tt d}(j',o)\geq n,\\
\big[\vartheta ({\tt d}(j,o)-\ovr,n-\ovr)\\
\quad -\vartheta ({\tt d}(j',o)-\ovr,n-\ovr)\big]^2,\\
\qquad\hbox{ if }\;\ovr<{\tt d}(j,o),{\tt d}(j',o)\leq n,\\
\vartheta ({\tt d}(j,o)-\ovr,n-\ovr)^2,\\
\qquad\hbox{ if }\;\ovr<{\tt d}(j,o)\leq n,
{\tt d}(j',o)\in ]\ovr,n[,\\
\vartheta ({\tt d}(j',o)-\ovr,n-\ovr)^2,\\
\qquad\hbox{ if }\;\ovr<{\tt d}(j',o)\leq n,
{\tt d}(j,o)\in ]\ovr,n[.\end{cases}\eqno (4.2.10)$$

By using convexity of the function exp and Eqn (4.2.9), $\forall$ $a>1$,
$$\begin{array}{l}\diy\frac{a}{2}\exp\,\Big[- 
h^{\Lam (n)}\Big({\bttg}_{\Lam (n)}\big(
\obom^0\vee
\bom_{\Lam (n)\setminus\Lam^0}\big)|
\bom_{\Gam\setminus\Lam (n)}\Big)\Big]\\
\quad +\diy\frac{a}{2}\exp\,\Big[- 
h^{\Lam (n)}\Big({\bttg}_{\Lam (n)}^{-1}\big(
\obom^0\vee
\bom_{\Lam (n)\setminus\Lam^0}\big)| 
\bom_{\Gam\setminus\Lam (n)}\Big)\Big]\\
\;\;\diy\geq a\exp\;\bigg[-\frac{1}{2}
h^{\Lam (n)}\Big({\bttg}_{\Lam (n)}\big(
\obom^0\vee
\bom_{\Lam (n)\setminus\Lam^0}\big), 
\bom_{\Gam\setminus\Lam (n)}\Big)\\
\qquad -\diy\frac{1}{2}
h^{\Lam (n)}\Big({\bttg}_{\Lam (n)}^{-1}\big(
\obom^0\vee
\bom_{\Lam (n)\setminus\Lam^0}\big)|
\bom_{\Gam\setminus\Lam (n)}\Big)\bigg]\\
\;\;\geq a\exp\,\Big[
-h^{\Lam (n)}\Big(
\obom^0\vee\bom_{\Lam (n)\setminus\Lam^0}|
\bom_{\Gam\setminus\Lam (n)}\Big)\Big]
e^{-C\Psi /2}\end{array}\eqno (4.2.11)$$
where
$$\Psi =\Psi (n,{\ttg})=|\utheta|^2\sum_{(j,j')\in\Lam (n)\times\Gam}
J({\ttd}(j,j'))\left|\upsilon (n,j)-\upsilon (n,j')
\right|^2.\eqno (4.2.12)$$

The next remark is that 
$$\begin{array}{l}\Psi\leq 3|\utheta |^2
\sum\limits_{(j,j')\in\Lam (n)\times\Gam}
{\mathbf 1}\big({\tt d}(j,o)\leq {\tt d}(j',o)
\big)J_{j,j'}\\ \;\\
\qquad\times\Big[\vartheta ({\tt d}(j,o)-\ovr,n-\ovr) -
\vartheta ({\tt d}(j',o)-\ovr,n-\ovr)\Big]^2\end{array}\eqno (4.2.13)
$$
where, with the help of the triangle inequality, for all 
$j,j': {\tt d}(j,o)\leq {\tt d}(j',o)$
$$\begin{array}{r}
0\leq\vartheta ({\tt d}(j,o)-\ovr,n-\ovr) -
\vartheta ({\tt d}(j',o)-\ovr,n-\ovr)\qquad{}\\
\leq {\tt d}(j,j')\diy\frac{z({\tt d}(j,o)-\ovr)}{Q(n-\ovr)}.
\end{array}\eqno (4.2.14)$$
This yields
$$\begin{array}{l}
\Psi\leq\diy\frac{3||\utheta||^2}{Q(n-\ovr)^2}\sum\limits_{(j,j')\in\Lam (n)\times\Gam}
J({\tt d}(j,j')){\tt d}(j,j')^2z({\tt d}(j,o)-\ovr)^2\\
\quad\leq\diy\frac{3||\utheta||^2}{Q(n-\ovr)^2}
\Big[\sup_{j\in\Gam}\sum\limits_{j'
\in\Gam}J({\tt d}(j,j')){\tt d}(j,j')^2 \Big]\sum\limits_{j\in\Lam_{n+r_0}}
z({\tt d}(j,o)-\ovr)^2.\end{array}
$$

In view of (1.3.4) it remains to estimate the sum $\sum\limits_{j\in\Lam_{n+r_0}}
z({\tt d}(j,o)-\ovr)^2$. To this end, observe that 
$uz(u)<1$ when $u\in (3,\infty )$. The next remark is that the
number of sites in the sphere $\Sigma_n$ grows linearly with $n$.
Consequently, 
$$\begin{array}{l}\sum\limits_{j\in\Lam(n+r_0)}
z({\tt d}(j,o)-\ovr)^2=\sum\limits_{1\leq k\leq n+r_0}z(k-\ovr)
\sum\limits_{j\in\Sigma_k}z(k-\ovr)\\
\qquad\qquad\leq C_0\sum\limits_{1\leq k\leq n+r_0}z(k-\ovr)
\leq C_1Q(n+r_0-\ovr)\end{array}$$
and
$$\Psi\leq\frac{C}{Q(n-\ovr)}\to\infty,\;\hbox{ as }\;n\to\infty.$$

Therefore, given $a>1$ for $n$ large enough,
the term $ae^{-C\Psi /2}$ in the RHS of (4.2.11)
becomes $>1$. Hence,
$$\begin{array}{l}\diy\frac{a}{2}\exp\,\Big[- 
h^{\Lam (n)}\Big({\bttg}_{\Lam (n)}\big(
\obom^0\vee\bom_{\Lam (n)\setminus\Lam^0}\big)
|{\mbox{\boldmath$\oomega$}}_{\Gam\setminus\Lam (n)}\Big)\Big]\\
\quad +\diy\frac{a}{2}\exp\,\Big[- 
h^{\Lam (n)}\Big({\bttg}_{\Lam (n)}^{-1}\big(
\obom^0\vee
\bom_{\Lam (n)\setminus\Lam^0}\big)|
\bom_{\Gam\setminus\Lam (n)}\Big)\Big]\\
\qquad\qquad\qquad
\;\;\geq\exp\,\Big[
-h^{\Lam (n)}\Big(\obom^0\vee
\bom_{\Lam (n)\setminus\Lam^0}|
\bom_{\Gam\setminus\Lam (n)}\Big)\Big]\,.
\end{array}\eqno (4.2.15)$$

Eqn (4.2.15) implies that  the conditional RDMF 
$$\begin{array}{r}
q^{\Lam^0|\Gam\setminus\Lam (n)}
(\obom^0|\bom_{\Gam\setminus\Lam (n)})=
\diy\int_{W_{\Lam (n)\setminus\Lam^0}}{\rd}\nu_{\Lam (n)\setminus\Lam^0}
(\bom_{\Lam (n)\setminus\Lam^0}) 
\quad\qquad{}\\
\diy\times\frac{\exp\big[-h^{\Lam^0}
(\obom^0\vee\bom_{\Lam (n)\setminus\Lam^0}|\bom_{\Gam\setminus\Lam^0})
\big]}{\Xi_{\Lam (n)}(
\bom_{\Gam\setminus\Lam (n)})},\end{array}
\eqno (4.2.16)$$
obeys
$$\begin{array}{r}
\lim\limits_{n\to\infty}\Big[q^{\Lam^0|\Gam\setminus\Lam (n)}
({\ttg}\obom^0|\bom_{\Gam\setminus\Lam (n)})
+q^{\Lam^0|\Gam\setminus\Lam (n)}
({\ttg}^{-1}\obom^0|\bom_{\Gam\setminus\Lam (n)})\Big]\qquad{}\\
\geq 2\lim\limits_{n\to\infty}
q^{\Lam^0|\Gam\setminus\Lam (n)}
(\obom^0|\bom_{\Gam\setminus\Lam (n)})\end{array}\eqno (4.2.17)$$
uniformly in boundary condition $\bom_{\Gam\setminus\Lam (n)}$. 
Integrating (4.2.17) in\\ ${\rd}
\mu_{\Gam}^{\Gam\setminus\Lam (n)}
(\bom_{\Gam\setminus\Lam (n)})$ 
yields (4.2.3). $\quad\Box$
\vskip .5 truecm 

{\bf 4.3.  Proof of Theorem 3.4 and Lemma 1.1.} In Theorem 3.4 
our argument follows the idea proposed in \cite{ISV}. 
On the lattice ${\bbZ}^2$ consider the squares 
$$\Lam (n)=\{j=(j_1,j_2):\;\max\,\big[|j_1|,|j_2|\big]\leq n\}, n=1,2,\ldots$$
The outer boundary of $\Lam (n)$ is the set 
$$\Sigma (n+1)=\{j=(
j_1,j_2):\;\max\,\big[|j_1|,|j_2|\big]=n+1\}.$$ 
Fix a point $x^*\in S^1$, 
a value $\beta\in (0,\infty )$ and 
consider a state $\vphi^*=\vphi (\mu^*)$ 
induced by measure $\mu^*=\mu^{(x^*)}\in{\fG}(\beta )$ which is a limiting 
point for the family of measures $\mu^*_n=\mu^{(x^*)}_{\beta,\Lam (n)}$ 
(cf. (2.3.10)) as $n\to\infty$. Here
$\mu^*_n$ stands for the probability distribution on 
$W_{\Lam (n)}$ with the `cooled' boundary condition
$$\bom^*_{\Sigma (n+1)}
=\{\om^*_j,\;j\in\Sigma (n+1)\}\eqno (4.3.1)$$
where
$$\om^*_j(\tau )\equiv x^*,\;\;0\leq\tau\leq\beta .\eqno (4.3.2)$$

Without loss of generality, assume that $\mu^*=\lim\limits_{n\to\infty}
\mu^*_n$. To simplify the notation, let us also omit the subscript $^*$,
writing $\mu =\mu^*$. 
If state $\vphi =\vphi^*$ is not $S^1$-invariant, we are done. 
Otherwise, suppose that $\vphi$ is $S^1$-invariant. Then choose 
an arc $\alpha =(x^*-1/200,x^*+1/200)$
of length $1/100$ around the point $x^*$ and let $\Pi_0(\alpha )$ 
be the orthoprojection on
the subspace in $\cH_0\simeq\cH$ formed by functions supported
by arc $\alpha$. Then 
$$\vphi (\Pi_0(\alpha ))=\int_{W_{\{0\}}}
{\rd}\mu^{\{0\}}(\omega_0){\mathbf 1}(x_0(\om_0)\in\alpha )
=\frac{1}{100}.\eqno (4.3.3)$$
(The lower/upper scripts $0$ and $\{0\}$ 
indicate that we take $\Lam^0=\{0\}$,
i.e., consider spins attached to lattice site $0\in\bbZ^2$.) 

Hence, for $n$ large enough, the conditional distribution\\
${\rd}\mu^{(\{0\}|\Sigma (n+1))}(\omega_0|
\bom^*_{\Sigma (n+1)})$ for $\omega_0\in W_{\{0\}}$, given
boundary condition $\bom^*_{\Sigma (n+1)}$, satisfies:
$$\int_{W_{\{0\}}}
{\rd}\mu^{(\{0\}|\Sigma (n+1))}(\omega_0|
\bom^*_{\Sigma (n+1)}){\mathbf 1}(x_0(\om_0)\in\alpha )
<\frac{1}{99}.\eqno (4.3.4)$$

Next, given $\eta\in (0,1]$, consider a family of points 
$${\wtx}_{j,\eta}
=x^*+j_1\eta\theta\;{\rm{mod}}\;1,\;j=(j_1,j_2)\in{\bbZ}^2,$$ 
and the family of the
corresponding cooled loops $\{{\wt\omega}_{j,\eta}\}$:  
$${\wt\omega}_{j,\eta}(\tau )\equiv{\wtx}_{j,\eta},
\;\;0\leq\tau\leq\beta .\eqno (4.3.5)$$ 
Further, consider the loop configuration ${\wt\bom}_{\Sigma (n+1),\eta}
=\{{\wt\omega}_{j,\eta},\;j\in\Sigma (n+1)\}$ over the boundary $\Sigma (n+1)$
formed by loops ${\wt\omega}_{j,\eta}$. For $\eta =1$, the only 
configuration over $\Lam (n)$ compatible with the boundary condition  
${\wt\bom}_{\Sigma (n+1),\eta}$ is the one where all loops
coincide with ${\wt\omega}_{\eta ,j}$: for any other choice of
the configuration, the energy $h^{\Lam (n)|\Sigma (n+1)}$ is equal to 
$+\infty$.
By continuity, for each $n$, there exists ${\wt\eta}(n)\in (0,1)$ such that
the probability measures 
${\rd}\mu^{(\{0\}|\Sigma (n+1))}(\omega_0|
{\wt\bom}_{\Sigma (n+1),{\wt\eta}(n)})$ 
conditional on ${\wt\bom}_{\Sigma_{n+1},{\wt\eta}(n)})$ satisfy  
$$\int_{W_{\{0\}}}
{\rd}\mu^{(\{0\}|\Sigma (n+1))}(\omega (0)|
{\wt\bom}_{\Sigma (n+1),{\wt\eta}(n)}){\mathbf 1}
(x_0(\om_0)\in\alpha )=\frac{2}{3},\eqno (4.3.6)$$
and ${\wt\eta}(n)$ is uniformly separated from $0$ and $1$. Any limiting 
point ${\wt\mu}$ of the sequence of conditional measures
$\mu^{(\Lam (n)|\Sigma (n+1))}(\,.\,|
{\wt\bom}_{\Sigma (n+1),{\wt\eta}(n)})$, $n\to\infty$, yields 
$$\int_{W_{\{0\}}}
{\rd}{\wt\mu}^{(\{0\}}(\omega_0){\mathbf 1}
(x_0(\om_0)\in\alpha )=\frac{2}{3},\eqno (4.3.7)$$
and the induced state ${\wt\vphi}$ gives  
${\wt\vphi}(\Pi_{\alpha}(0))=2/3$. It means that neither
${\wt\mu}$ nor ${\wt\vphi}$ are $S^1$-invariant. $\quad\Box$
\vskip .5 truemm

{\it Proof of Corollary} 3.5. Eqn (4.3.7) guarantees that
$\forall$ $\beta >0$ there exists a non-$S^1$-invariant
measure ${\wt\mu}_\beta\in\fG(\beta )$. Passing to a limiting point
as $\beta\to\infty$ yields a ground-state measure ${\wt\psi}$ with
the property that 
$$\int_{W_{\{0\}}}
{\rd}{\wt\psi}^{\{0\}}(\omega_0){\mathbf 1}
(x_0(\om_0)\in\alpha )=\frac{2}{3},\eqno (4.3.8)$$
again contradicting $S^1$-invariance. $\quad\Box$
\vskip .5 truecm

{\textit{Proof of Lemma 1.1.}} Let $\lambda_1\geq\lambda_2\geq\ldots$ be the
sequence of the eigenvalues of operator $R$ and $e_i(x)$, $
i=1,2,\ldots$ be the corresponding eigenvectors. As follows
from (1.4.11),
$$\lim_{n\to\infty}\sum_{i,j}\Big(\big\langle R_ne_i,e_j
\big\rangle -\lambda_i\delta_{ij}\Big)^2=0.$$

We want to show that the sequence $\{R_n\}$ converges to $R$ in
the Banach space ${\mathfrak C}$ of the trace-class operators
in $L_2(M,\nu )$ with the norm $\|\,\cdot\,\|_1$.
We are going to use a natural basis in ${\mathfrak C}$
formed by the system of rank one `matrix units' $E_{ij}=
|e_i\rangle\langle e_j|$. Set:
$$\begin{array}{c}R^{(i_0)}_n=\sum\limits_{1\leq i,j<i_0}
\big\langle R_ne_i,e_j\big\rangle E_{ij},\;\;
{\overline R}^{(i_0)}_n=\sum\limits_{i,j\geq i_0}
\big\langle R_ne_i,e_j\big\rangle E_{ij},\\
{\widetilde R}^{(i_0)}_n=\sum\limits_{1\leq i<i_0}
\sum\limits_{j\geq i_0}
\big\langle R_ne_i,e_j\big\rangle E_{ij}.\end{array}$$
Next, set
$$R^{(i_0)}=\sum_{1\leq i<i_0}\lambda_iE_{ii}.$$

Clearly, $R^{(i_0)}_n$ and ${\overline R}^{(i_0)}_n$ are
positive-definite operators. Furthermore,
$$\big\|R^{(i_0)}_n\big\|_1+\big\|{\overline R}^{(i_0)}_n
\big\|_1=1$$
and
$$R_n=R^{(i_0)}_n+{\overline R}^{(i_0)}_n
+{\wt R}^{(i_0)}_n+\left({\wt R}^{(i_0)}_n\right)^*.$$
Take an arbitrary $\epsilon >0$ and choose $i_0=i_0(\epsilon )$
and $n_0=n_0(\epsilon )$ such that
$$\sum_{i\geq i_0}\lambda_i<\frac{\epsilon}{8}$$
and for $n\geq n_0$
$$\left|\|R_n^{(i_0)}-R^{(i_0)}\right\|_1
<\frac{\epsilon}{8},\;\;\sum_{i\neq j}
\left(\left\langle R_ne_i,e_j\right\rangle\right)^2
<\frac{\epsilon^2}{{\sqrt 2}i_0^2}.$$
Then for $n\geq n_0$,
$$\begin{array}{cl}
\left\|R-R_n\right\|_1&
\leq\left\|R-R^{(i_0)}\right\|_1
+\left\|R^{(i_0)}-R_n^{(i_0)}\right\|_1
+\left\|R_n^{(i_0)}-R_n\right\|_1\\
\;&\leq\epsilon/8 +\epsilon/8 +
\left\|R_n-R_n^{(i_0)}\right\|_1.\end{array}$$

It remains to estimate the term $\left\|R_n-R_n^{(i_0)}\right\|_1$.
To this end we write:
$$\begin{array}{cl}
\left\|R_n-R_n^{(i_0)}\right\|_1&\leq
\left\|{\overline R}_n^{(i_0)}\right\|_1
+2\left\|{\widetilde R}_n^{(i_0)}\right\|_1
=1-\left\|R_n^{(i_0)}\right\|_1
+2\left\|{\widetilde R}_n^{(i_0)}\right\|_1\\
\;&\leq 1-\left\|R^{(i_0)}\right\|_1+\epsilon/8
+2\left\|{\widetilde R}_n^{(i_0)}\right\|_1\leq 1-\epsilon/4
+2\left\|{\widetilde R}_n^{(i_0)}\right\|_1.\end{array}
$$
Finally,
$$\left\|{\widetilde R}_n^{(i_0)}\right\|_1\leq
\sum_{1\leq i<i_0}
\left[\sum_{j\geq i_0}\langle R_ne_i,e_j\rangle^2\right]^{1/2}
<i_0\left[\sum_{i\neq j}\langle R_ne_i,e_j\rangle^2\right]^{1/2}$$
which is $<\epsilon/2$. This completes the proof of Lemma 1.1. 
$\quad\Box$
\vskip 1 truecm

{\bf Acknowledgement.}   
This work has been conducted under Grant \\
2011/20133-0 provided by 
the FAPESP, Grant 2011.5.764.45.0 provided by The Reitoria of the 
Universidade de S\~{a}o Paulo and Grant 2012/04372-7
provided by the FAPESP. The authors
express their gratitude to NUMEC and IME, Universidade de S\~{a}o Paulo,
Brazil, for the warm hospitality. The authors thank the referees for remarks
and suggestions.

\vfill\eject

\end{document}